\documentclass[a4paper, reqno]{amsart}
\pdfoutput=1
\usepackage[english]{babel}
\usepackage{amssymb, mathtools, braket, upref}
\usepackage{tikz, pgfplots}
\pgfplotsset{compat=newest}
\usepackage[shortlabels]{enumitem}
\usepackage{xcolor}
\usepackage[colorlinks=true]{hyperref}
\hypersetup{urlcolor=blue, citecolor=red}

\title[Heat equation with non-local Robin conditions]{Non-Positivity of the heat equation with non-local Robin boundary conditions}

\author[J. Gl\"{u}ck]{Jochen Gl\"{u}ck}
\address{Jochen Gl\"{u}ck, Bergische Universit\"{a}t Wuppertal, Gaussstr. 20, 42119, Wuppertal, Germany}
\email{glueck@uni-wuppertal.de}

\author[J. Mui]{Jonathan Mui}
\address{Jonathan Mui, Bergische Universit\"{a}t Wuppertal, Gaussstr. 20, 42119, Wuppertal, Germany}
\email{jomui@uni-wuppertal.de}


\date{\today}

\numberwithin{equation}{section}

\theoremstyle{plain}
\newtheorem{theorem}{Theorem}[section]
\newtheorem{proposition}[theorem]{Proposition}
\newtheorem{lemma}[theorem]{Lemma}
\newtheorem{corollary}[theorem]{Corollary}

\theoremstyle{definition}

\newtheorem{example}[theorem]{Example}
\newtheorem{examples}[theorem]{Examples}
\newtheorem{setting}[theorem]{Setting}

\theoremstyle{remark}
\newtheorem{remark}[theorem]{Remark}

\DeclareMathOperator{\dv}{div}
\DeclareMathOperator{\dom}{dom}

\DeclareMathOperator{\id}{id}

\newcommand{\mf}{\mathfrak}
\newcommand{\rs}{\mathrm{s}}
\newcommand{\NN}{\mathbb{N}}
\newcommand{\ZZ}{\mathbb{Z}}
\newcommand{\RR}{\mathbb{R}}
\newcommand{\CC}{\mathbb{C}}
\newcommand{\modulus}[1]{\left\lvert #1 \right\rvert}
\newcommand{\argument}{\mathord{\,\cdot\,}}

\renewcommand{\Re}{\operatorname{Re}}
\renewcommand{\Im}{\operatorname{Im}}

\newcommand{\norm}[1]{\left\lVert #1 \right\rVert}
\newcommand{\one}{\mathbf{1}}

\begin{document}
	
\begin{abstract}
    We study heat equations $\partial_t u - \dv(A\nabla u) = 0$ on bounded Lipschitz domains $\Omega$, where $-\dv(A\nabla\,\cdot\,)$ is a second-order uniformly elliptic operator with generalised Robin boundary conditions. 
    These boundary conditions are formally given by $\nu\cdot A\nabla u + Bu=0$, 
    where $B\in\mathcal{L}(L^2(\partial\Omega))$ is a general operator. 
    In contrast to large parts of the literature on non-local Robin boundary conditions, we also allow for operators $B$ that destroy the positivity preserving property of the solution semigroup. 
    Nevertheless, we obtain ultracontractivity of the semigroup under quite mild assumptions on $B$. 
    For a certain class of operators $B$ we demonstrate that the semigroup is in fact eventually positive rather than positivity preserving.
\end{abstract}

\maketitle
    
\section{Introduction}

\subsection{Main results and outline of paper}

We begin with an informal description of the topic of the article. Let $\Omega \subseteq \RR^d$ be a bounded Lipschitz domain, and let $A \colon \Omega \to \RR^{d\times d}$ be a matrix of coefficients consisting of real-valued, bounded measurable functions.
We study solutions $u=u(t,x)$ to the uniformly elliptic second-order parabolic equation
\begin{equation}
    \label{eq:pde}
    \frac{\partial u}{\partial t}  - \dv(A\nabla u) = 0 \qquad \text{for } (t,x)\in (0,\infty)\times\Omega
\end{equation}
with a generalised Robin boundary condition formally given by
\begin{equation}
    \label{eq:Robin-boundary-cond}
	\nu\cdot A\nabla u + Bu = 0 \qquad\text{on } \partial\Omega,
\end{equation}
where $\nu$ is the outer unit normal on $\partial\Omega$ and $B$ is a bounded linear operator on $L^2(\partial\Omega)$.
The classical Robin boundary conditions are recovered by taking $B$ to be a multiplication operator given by a real-valued function $\beta\in L^\infty(\partial\Omega)$. On the other hand, the general form of the boundary operator we consider allows the possibility of \emph{non-local} boundary conditions, for instance if $B$ is an integral operator
\begin{equation*}
	(Bf)(x) \coloneqq \int_{\partial\Omega} k(x,y)f(y)\,dy \qquad\forall\, f\in L^2(\partial\Omega)
\end{equation*}
induced by a measurable function $k \colon \partial\Omega\times\partial\Omega \to \CC$. The differential operator $- \dv(A\nabla \argument)$ in~\eqref{eq:pde} that is subject to the generalised Robin boundary condition~\eqref{eq:Robin-boundary-cond} is rigorously defined as a closed operator $L_B$ on $L^2(\Omega)$ associated to a sesquilinear form $\mf{a}_B$; see Setting~\ref{setting:main} in Section~\ref{sec:forms} below.

Many qualitative properties of solutions to~\eqref{eq:pde} and~\eqref{eq:Robin-boundary-cond} are most easily formulated by considering the semigroup $(e^{-tL_B})_{t\ge 0}$ generated by $-L_B$. 
The main theme in our investigation is the lack of positivity of the semigroup. 
For classical boundary conditions (Dirichlet, Neumann, Robin, mixed) --- which, to emphasise, are local --- it is well-known that a \emph{positivity preserving property} holds for the evolution equation~\eqref{eq:pde}. This means that given an initial function $u_0 \ge 0$, the corresponding solution $u$ satisfies $u(t,x)\ge 0$ for all $x\in\Omega$ and for all $t>0$. 
For non-local Robin boundary conditions, this does not hold in general, and in fact the positivity preserving property can be easily characterised in terms of $B$, see Proposition~\ref{prop:positivity} below. 
In the positivity preserving case, there is a substantial body of work on non-local Robin boundary conditions, even for operators $B$ that are unbounded on the boundary space $L^2(\partial \Omega;\CC)$; see Section~\ref{sec:earlier} for details.

In contrast, we are mostly interested in the case where positivity is not preserved. 
We focus on two questions that arise in this situation. 
On the one hand, we study whether one still has \emph{ultracontractivity} of the semigroup $(e^{-tL_B})_{t\ge 0}$, which is the property that for each $t>0$, the operator $e^{-tL_B}$ maps $L^2$ into $L^\infty$. 
This property is commonly shown by combining a Sobolev embedding theorem with an interpolation result which requires the semigroup to be bounded for small times on the spaces $L^1$ and $L^\infty$. Without positivity, this boundedness is not straightforward to obtain.
On the other hand, given that the semigroup $(e^{-tL_B})_{t\ge 0}$ will not preserve positivity in general, we give sufficient conditions for the weaker property of \emph{uniform eventual positivity} that is explained after Theorem~\ref{thm:intro} and in Section~\ref{sec:epos}. 

To give the reader an overview of the main results, we state the following theorem, which is a simplified combination of Theorems~\ref{thm:ultra} and~\ref{thm:eventual-pos}. 
Inequalities between $L^2$-functions are meant almost everywhere; see Section~\ref{subsection:real-and-positive} for details.

\begin{theorem}
    \label{thm:intro}
    Let $\Omega\subseteq\RR^d$ be a bounded Lipschitz domain, and let $B$ be a bounded and self-adjoint linear operator on $L^2(\partial\Omega)$ that leaves $L^\infty(\partial\Omega)$ invariant and maps real-valued functions to real-valued functions.

    The semigroup $(e^{-tL_B})_{t\ge 0}$ enjoys the following properties:
    
    \begin{enumerate}[\upshape(i)]
	    \item\label{thm:intro:it:ultracon} 
        The semigroup is \emph{ultracontractive}, i.e.\ 
        \[ 
            e^{-tL_B}\big(L^2(\Omega)\big) \subseteq L^\infty(\Omega) \quad\forall\, t>0.
        \]
    	
        \item\label{thm:intro:it:unif-evpos} 
        If $B$ is positive semi-definite and $B\one_{\partial\Omega}=0$, then there exist $t_0\ge 0$ and $\delta>0$ such that
    	\begin{equation*}
    		e^{-tL_B}f \ge \delta\left(\int_\Omega f\,dx\right)\one \qquad\forall\, t\ge t_0
    	\end{equation*}
        for every $0\le f\in L^2(\Omega)$.
    \end{enumerate}
\end{theorem}

Here, $\one$ denotes the constant function with value $1$ on $\Omega$ and $\one_{\partial \Omega}$ denotes the constant function with value $1$ on $\partial \Omega$.

The conclusion of part~\ref{thm:intro:it:unif-evpos} of the theorem implies, in particular, that $e^{-tL_B}f \ge 0$ for all $t\ge t_0$ and all $0 \le f \in L^2(\Omega)$. 
This property is called \emph{uniform eventual positivity} of the semigroup $(e^{-tL_B})_{t\ge 0}$; see Section~\ref{sec:epos} for more information.

In Section~\ref{sec:fun}, we recall some background material, define the differential operator $L_B$ via a sesquilinear form, and collect some functional analytic properties of $L_B$ and its associated semigroup $(e^{-tL_B})_{t\ge 0}$ on $L^2(\Omega)$. Section~\ref{sec:ultra} is devoted to ultracontractivity of the semigroup (Theorem~\ref{thm:ultra}). For our results on eventual positivity, we require some spectral conditions on $B$, which turn out to be directly related to spectral properties of the differential operator $L_B$. In the present paper, we consider two simple conditions: firstly, the case $B\one_{\partial\Omega} = 0$ (as in Theorem~\ref{thm:intro}\ref{thm:intro:it:unif-evpos} above) will be discussed in Section~\ref{sec:spectrum-epos}, and the condition $\braket{B\one_{\partial\Omega}, \one_{\partial\Omega}}<0$ will be treated in Section~\ref{sec:symm}. In the latter section, we focus the analysis on the special case that $\Omega$ is a ball.

\subsection{Earlier work on non-local Robin problems}
\label{sec:earlier}

Non-local Robin boundary conditions appear in the literature as far back as the 1950's, due to Feller in his seminal work on Markov diffusions in one dimension~\cite{F54}. Thanks to the probabilistic connection, positive semigroups arise very naturally in this context. 
Closer to the current day, specific examples of non-local Robin conditions have appeared in the study of Schr\"{o}dinger operators~\cite{Schroe88}, a model of Bose condensation~\cite{Schroe89}, a reaction-diffusion equation~\cite{GuMe97}, and in a model of a thermostat~\cite{GuMe00}. The latter paper demonstrated explicitly that the associated semigroup could fail to be positive. 
One of the earliest treatments of non-local Robin conditions in a general functional analytic framework appears in the work of Gesztesy and Mitrea~\cite{GeMi08}. The development of the abstract theory continued in collaboration with other co-authors in~\cite{GeMi09, GeMiNi14, GeMiNiOu}, with a particular focus on sesquilinear forms, positive semigroups and Gaussian estimates. 

Non-local Robin conditions of a different type, closer to Feller's original inspiration, appear in the work of Arendt et al.~\cite{AKK}. Using modern developments in semigroup theory, the authors were able to deduce smoothing properties (the strong Feller property and holomorphy), contractivity, and analyse the asymptotic behaviour of the semigroup. Further variations on the boundary operators and even extensions to nonlinear equations may be found in~\cite{Ait21, Velez-Warma10, Velez11}.
	
Except for~\cite{GuMe00, GeMiNiOu}, the works mentioned above all feature Robin boundary conditions that produce positive semigroups. 
The more subtle property of eventual positivity was first analysed in the papers~\cite{DGK2, DG18}, within a general theory of eventually positive semigroups on Banach lattices. The specific models of the thermostat and Bose condensation (\cite{GuMe00} and~\cite{Schroe89} respectively) are revisited in~\cite[Section 6]{DGK2}, where conditions on the boundary operator $B$ are given such that the associated semigroup is non-positive but eventually positive.

\subsection{A taste of eventual positivity}
\label{sec:epos}

A $C_0$-semigroup $(e^{-tC})_{t \ge 0}$ on an $L^p$-space (or, more generally, on a Banach lattice) is called \emph{individually eventually positive} if for each $f \ge 0$ there exists a time $t_0 \ge 0$ such that $e^{-tC}f \ge 0$ for all $t \ge t_0$. 
If $t_0$ can be chosen to be independent of $f$, then the semigroup is called \emph{uniformly eventually positive}.  
In the situation of Theorem~\ref{thm:intro}\ref{thm:intro:it:unif-evpos}, the theorem implies that $(e^{-tL_B})_{t \ge 0}$ is uniformly eventually positive. 

There is a variety of stronger notions of eventual positivity.
For instance, in the terminology of \cite{DGK2}, the conclusion of Theorem~\ref{thm:intro}\ref{thm:intro:it:unif-evpos} even implies that the semigroup is \emph{uniformly eventually strongly positive with respect to $\one$}. 
However, the details of this nomenclature are quite technical and mainly useful to develop and structure the general theory of eventual positivity, which is not the goal of the present paper. 
In order not to overload the paper with terminology, we thus refrain from discussing all these notions in detail. 
Instead, we simply state the lower estimates in each result explicitly and restrict ourselves to mentioning that those estimates imply, in particular, that the semigroup is uniformly eventually positive.

We briefly comment on the history of eventual positivity.
The phenomenon has been known for quite some time in finite dimensions~\cite{NT} and in some concrete examples of partial differential equations, e.g.\ in fourth-order parabolic equations~\cite{FGG}. In~\cite{Dan-DN}, Daners investigated eventual positivity for the semigroup generated by the Dirichlet-to-Neumann operator on the unit disk. 
This case study motivated the development of the general theory of eventually positive semigroups using abstract techniques from Banach lattice and operator theory, in collaboration with Kennedy and the first-named of the present authors, in the articles~\cite{DGK1, DGK2, DG18}. Since the publication of these works, the theory on eventual positivity has branched off in various directions. The interested reader may consult the survey article~\cite{G22} for a `bird's-eye view' of the subject of eventual positivity and more references to recent developments. In particular, we mention that the functional analytic approach to eventual positivity has proved to be especially useful for evolution equations with higher-order differential operators, see e.g.~\cite{AGRT} and~\cite{DKP}.

One of the core ingredients in the abstract study of eventual positivity is the spectral theory of positive operators, motivated by two celebrated results: the Perron-Frobenius theorems in finite dimensions and the Krein-Rutman theorem in infinite dimensions. 
The existence of a positive leading eigenfunction of the differential operator and certain spectral considerations are crucial in order to apply the results of the theory. 
Another essential ingredient for the theory in infinite dimensions is a smoothing condition on the semigroup, which often translates to ultracontractivity in PDE applications. 
The present article will demonstrate both of the core ingredients in action. 
This feature of the general theory of eventually positive semigroups explains the specific spectral assumptions on the boundary operator $B$ that we will consider in Sections~\ref{sec:spectrum-epos} and~\ref{sec:symm}, and also the need for an ultracontractivity result. 
Such a result is implicitly contained in \cite[Theorem~3.6]{GeMiNiOu}, since Gaussian estimates for the semigroup are proved there. However, the assumptions in this theorem are considerably different from ours; see the discussion after Proposition~\ref{prop:local-by-dom} for details.

Finally, one could wonder what happens if one or both of the core ingredients mentioned above are not available. This can happen, for instance, if one is interested in higher-order evolution equations on unbounded domains (spectral conditions fail) or non-smooth domains (smoothing condition fails). In such cases, the development of a general theory is still very much open; however, recent work of Arora~\cite{Ar22} and the second author of the present paper~\cite{Mui} represents some steps in this direction.

\section{Setting the stage}
\label{sec:fun}

In this preliminary section, we collect some basic facts that will be essential to our analysis, and also fix notations and terminology.
In the sequel, we consider complex-valued functions unless otherwise stated. If $X$ is a Banach space, we write $\mathcal{L}(X)$ for the space of bounded linear operators on $X$.

\subsection{Lipschitz domains}
\label{sec:Lipschitz}
In this article, a \emph{domain} $\Omega\subseteq\RR^d$ is a non-empty, connected, open set. We say that $\Omega$ is a \emph{Lipschitz domain} if its boundary $\partial\Omega$ is locally the graph of a Lipschitz function. Precisely, this means the following: for all $x_0\in\partial\Omega$, there exist $\delta,\varepsilon>0$, an orthogonal transformation $T \colon \RR^d\to \RR^d$ and a Lipschitz function $\varphi \colon \{x'\in\RR^{d-1}:\modulus{x'} < \delta\} \to \RR$
such that
\begin{equation*}
	U \coloneqq T^{-1}\big( \{(x', x_d)\in \RR^{d-1}\times\RR : \modulus{x'} < \delta, \modulus{x_d - \varphi(x')}< \varepsilon\} \big)
\end{equation*}
is a neighbourhood of $x_0$, and
\begin{equation*}
	\Omega\cap U = T^{-1}\Big( \{(x',x_d)\in \RR^{d-1}\times\RR : \modulus{x'} < \delta, 0 < x_d - \varphi(x') < \varepsilon\} \Big).
\end{equation*}

On a Lipschitz domain $\Omega \subseteq \RR^d$, we denote the $L^2(\Omega)$ norm by $\|\cdot\|_2$. 
The boundary $\partial\Omega$ is equipped with the $(d-1)$-dimensional Hausdorff measure $\sigma$, and the $L^2(\partial\Omega)$ norm is denoted by $\|\cdot\|_{2,\partial\Omega}$. 
The inner products on $L^2(\Omega)$ and $L^2(\partial\Omega)$ are respectively written as $\braket{\cdot,\cdot}_\Omega$ and $\braket{\cdot,\cdot}_{\partial\Omega}$.
Throughout, the constant function of value $1$ on $\Omega$ is denoted  by $\one$, 
and the constant function of value $1$ on $\partial \Omega$ is denoted by $\one_{\partial \Omega}$.

\subsection{Real and positive functions and operators} 
\label{subsection:real-and-positive}

Let $(M,\mu)$ be a $\sigma$-finite measure space and $p \in [1,\infty]$, and consider $E = L^p(M,\mu)$. 
Given $f,g\in E$, we write the inequality $f \le g$ to mean that $f, g$ are real-valued and $f(x) \le g(x)$ for $\mu$-a.e.\ $x\in M$. 
For a real-valued function $f \in E$ we define the \emph{positive part} $f^+$ by $f^+(x) \coloneqq \max\{f(x),0\}$ for $\mu$-a.e.\ $x \in M$, and the \emph{negative part} as $f^- \coloneqq (-f)^+$. 
A linear operator $T \colon E \to E$ is called \emph{positive} if $Tf \ge 0$ for all $0\le f\in E$.

Let $E_\RR$ denote the set of real-valued functions in $E$.
A closed linear operator $T:\dom(T)\subseteq E\to E$ is said to be \emph{real} if
\begin{equation*}
	\dom(T) = \dom(T)\cap E_\RR + i \dom(T)\cap E_\RR \quad\text{and}\quad T(\dom(T)\cap E_\RR)\subseteq E_\RR.
\end{equation*}
If $T \in \mathcal{L}(E)$, then the above conditions reduce to the single condition $T(E_\RR) \subseteq E_\RR$. Clearly, positive operators are real.

\subsection{Semigroup theory}
\label{sec:semigroups}

We assume that the reader is familiar with the theory of one-parameter operator semigroups. Recall that the \emph{spectral bound} of a closed operator $C \colon \dom(C)\subseteq X\to X$ on a Banach space $X$ is defined by
\begin{equation*}
	\rs(C) \coloneqq \sup\{\Re\lambda : \lambda\in\sigma(C)\} \in [-\infty,\infty],
\end{equation*}
where $\sigma(C) \subseteq \CC$ denotes the \emph{spectrum} of $C$.
If $-C$ is the generator of a strongly continuous semigroup $(e^{-tC})_{t\ge 0}$ on $X$, then the \emph{growth bound} of the semigroup is the quantity
\begin{equation*}
	\omega_0(-C) \coloneqq \inf\left\{ \omega\in\RR ~\colon \exists\, M \ge 1 \; \forall\,t\ge 0 \;\|e^{-tC}\|_{\mathcal{L}(X)} \le M e^{\omega t} \right\}.
\end{equation*}
It is a standard fact that the resolvent operator $R(\lambda,-C)\coloneqq(\lambda + C)^{-1}$ can be represented by the Laplace transform of the semigroup whenever $\Re\lambda>\omega_0(-C)$; namely
\begin{equation*}
	R(\lambda, -C)f = \int_0^\infty e^{-\lambda t} e^{-tC}f\,dt
\end{equation*}
for all $f\in X$ and $\lambda\in\CC$ with $\Re\lambda>\omega_0(-C)$, where the integral converges as a Bochner integral. We will use this fact in the proof of Proposition~\ref{prop:Linfty-bounded-pos-sgrp} below.

If $(M,\mu)$ is a $\sigma$-finite measure space and $E = L^p(M,\mu)$, $p\in [1,\infty]$, a $C_0$-semigroup $(e^{-tC})_{t\ge 0}$ on $E$ is called \emph{positive} if each operator $e^{-tC}$ is a positive operator on $E$. Similarly, we say that the semigroup $(e^{-tC})_{t\ge 0}$ is \emph{real} if each operator $e^{-tC}$ is real. Note that the semigroup $(e^{-tC})_{t\ge 0}$ is real if and only if $C$ is real. 

\subsection{Non-local Robin boundary conditions via forms}
\label{sec:forms}

Here we introduce the main setting which will be used throughout the rest of the paper, and provide the definition of the operator $L_B$ mentioned in the introduction. We will refer to this setting each time it is used.

\begin{setting}[Main Setting]
    \label{setting:main}
    Let $\Omega\subseteq\RR^d$ be a bounded Lipschitz domain. Let $A \colon \Omega \to \RR^{d\times d}$ be a matrix-valued map with the following properties:
    \begin{enumerate}[label=(A\arabic*)]
        \item\label{setting:main:itm:A-real} 
        The coefficients $a_{ij} \colon \Omega \to \RR$ of $A$ are bounded and measurable.
        
        \item\label{setting:main:itm:unif-elliptic} 
        $A$ is uniformly elliptic with lower bound $\alpha>0$, i.e.\ it satisfies
        \begin{equation*}
            \Re\left(\overline{\xi}^\top A(x)\xi\right) \ge \alpha \modulus{\xi}^2 \qquad \forall\, \xi\in\CC^d, x\in\Omega.
        \end{equation*}
    \end{enumerate}
    Let $B \in \mathcal{L}(L^2(\Omega))$.
    \begin{enumerate}[(a)]
        \item\label{setting:main:itm:form-b} 
        We define a bounded sesquilinear form $\mf{b} \colon L^2(\partial \Omega) \times L^2(\partial \Omega) \to \CC$ by
        \begin{equation*}
        	\mf{b}[f,g] \coloneqq \int_{\partial\Omega} (Bf)\overline{g}\,d\sigma = \braket{Bf,g}_{\partial\Omega}
        \end{equation*}
        for all $f,g\in L^2(\partial\Omega)$, 
        where $\sigma$ denotes the $(d-1)$-dimensional Hausdorff measure on $\partial\Omega$.

        \item\label{setting:main:itm:form-a}
        We define a sesquilinear form $\mf{a}_B$ on $L^2(\Omega)$ by
        \begin{align*}
        	\label{eq:a-form}
        	\dom(\mf{a}_B) &\coloneqq H^1(\Omega) \\
        	\mf{a}_B[u,v] &\coloneqq 
            \int_\Omega A\nabla u\cdot\overline{\nabla v}\,dx + \int_{\partial\Omega} (B\gamma(u))\overline{\gamma(v)}\,d\sigma \\
            & = 
            \int_\Omega A\nabla u\cdot\overline{\nabla v}\,dx + \mf{b}[\gamma(u),\gamma(v)] 
        \end{align*}
        for all $u,v\in H^1(\Omega)$, where $\gamma\in\mathcal{L}(H^1(\Omega), L^2(\partial\Omega))$ is the trace operator.

        \item\label{setting:main:itm:operator} 
        By $L_B \colon \dom(L_B) \subseteq L^2(\Omega) \to L^2(\Omega)$ we denote the operator associated to the form $\mf{a}_B$, i.e.\ the operator given by
        \begin{align*}
        	\dom(L_B) &= \left\{ u \in H^1(\Omega) ~\colon \exists\, f\in L^2(\Omega) \; \forall\,v\in H^1(\Omega) \; \, \mf{a}_B[u,v]=\braket{f,v}_{\Omega} \right\} \\
                L_B u &\coloneqq f,
        \end{align*}
        where the function $f \in L^2(\Omega)$ in the second line is that which occurs in the definition of $\dom(L_B)$ (it is unique since $\mf{a}_B$ is densely defined).

        \item\label{setting:main:itm:A-dependence}
        In the (rare) situations where we need to stress the dependence of the operator $L_B$ on the coefficient matrix $A$, we will use the notation $L_{A,B}$ instead of $L_B$, and refer explicitly to Setting~\ref{setting:main}\ref{setting:main:itm:A-dependence} when doing so.
    \end{enumerate}
\end{setting}

We collect some fundamental properties of the form $\mf{a}_B$ and its associated operator $L_B$ in the following proposition.

\begin{proposition}
\label{prop:basic}
	In Setting~\ref{setting:main}, the form $\mf{a}_B$ and its associated operator $L_B$ satisfy the following properties.
	\begin{enumerate}[\upshape(i)]
        \item\label{prop:basic:itm:H1-elliptic} 
        The form $\mf{a}_B$ is bounded on $H^1(\Omega)$ and $H^1$-elliptic; that is, there exist constants $c,\omega>0$ such that
		\begin{equation*}
			\Re\mf{a}_B[u,u] + \omega\|u\|^2_2 \ge c\|\nabla u\|^2_2 
            \qquad 
            \forall u \in H^1(\Omega)
            .
		\end{equation*}
        
		\item\label{prop:basic:itm:operator}
        
        The operator $L_B$ is densely defined, closed, has compact resolvent, and $-L_B$ generates an analytic $C_0$-semigroup $(e^{-tL_B})_{t\ge 0}$ on $L^2(\Omega)$.

        \item\label{prop:basic:itm:contractive} 
        If the operator $B+B^*$ on $L^2(\partial\Omega)$ is positive semi-definite, then the semigroup $(e^{-tL_B})_{t\ge 0}$ is contractive on $L^2(\Omega)$, i.e.\ $\|e^{-tL_B}\|_{L^2\to L^2} \le 1$ for all $t \ge 0$.
        
        \item\label{prop:basic:itm:adjoint}
        Using the notation $L_{A,B}$ for $L_B$ from Setting~\ref{setting:main}\ref{setting:main:itm:A-dependence}, one has $(L_{A,B})^* = L_{A^\top, B^*}$.

        \item\label{prop:basic:itm:self-adjoint}
        If the matrix $A(x)$ is symmetric for each $x \in \Omega$ and $B$ is self-adjoint, then the operator $L_{B}$ is self-adjoint.

        \item\label{prop:basic:itm:real}
        The operator $B \in \mathcal{L}(L^2(\partial\Omega))$ is real if and only if the semigroup $(e^{-tL_B})_{t\ge 0}$ is real.
	\end{enumerate}
\end{proposition}

\begin{proof}
	\ref{prop:basic:itm:H1-elliptic}
    The boundedness of the form $\mf{a}_B$ on $H^1(\Omega)$ follows easily from the boundedness of the coefficients of the matrix $A$ (assumption~\ref{setting:main:itm:unif-elliptic} in Setting~\ref{setting:main}) and from the continuity of the trace operator $\gamma \colon H^1(\Omega)\to L^2(\Omega)$.

    We now prove the $H^1$-ellipticity.
    By~\cite[Lemma 2.5]{GeMi08}, for every $\varepsilon>0$ there exists $C_\varepsilon>0$ such that
    \begin{equation}
        \|\gamma(u)\|^2_{2,\partial\Omega} \le \varepsilon \|\nabla u\|^2_2 + C_\varepsilon \|u\|^2_2 \qquad\forall\, u\in H^1(\Omega). 
    \end{equation}
    From this result, we obtain
	\begin{equation*}
		\modulus{\mf{b}[\gamma(u),\gamma(u)]} 
        \le 
        \|B\|_{\mathcal{L}(L^2(\partial\Omega))} (\varepsilon\|\nabla u\|^2_2 + C_\varepsilon\|u\|^2_2) 
        \quad \forall\, u\in H^1(\Omega).
	\end{equation*}
    If $\alpha$ is the ellipticity constant from~\ref{setting:main:itm:unif-elliptic}, it consequently holds that
	\begin{align*}
		\Re\mf{a}_B[u,u] &\ge \alpha \|\nabla u\|^2_2 + \Re\mf{b}[\gamma(u),\gamma(u)] \\
		&\ge \alpha \|\nabla u\|^2_2 - \|B\| (\varepsilon\|\nabla u\|^2_2 + C_\varepsilon\|u\|^2_2).
	\end{align*}
    For concreteness, we choose $\varepsilon\le \alpha (2\|B\|)^{-1}$ so that
    \begin{equation*}
        \Re\mf{a}_B[u,u] + \omega\|u\|^2_2 \ge \frac{\alpha}{2}\|\nabla u\|^2_2 \qquad\forall\, u\in H^1(\Omega)
    \end{equation*}
    with $\omega \coloneqq \|B\|C_\varepsilon$.

    \ref{prop:basic:itm:operator}
    The closedness of $L_B$ follows easily from boundedness of the form $\mf{a}_B$ on $H^1(\Omega)$ and the ellipticity proven above. The density of $\dom(L_B)$ in $L^2(\Omega)$ is standard, see e.g.~\cite[Proposition 1.22]{Ouh}.
    
    The compactness of the resolvent is a direct consequence of the embedding $H^1(\Omega) \hookrightarrow L^2(\Omega)$, which is compact since $\Omega$ is a bounded Lipschitz domain.

    The fact that $-L_B$ generates a strongly continuous analytic semigroup $(e^{-tL_B})_{t\ge 0}$ on $L^2(\Omega)$ follows by standard results on sesquilinear forms, e.g.~\cite[Theorem 1.52]{Ouh}.

    \ref{prop:basic:itm:contractive} 
	We first show that $\mf{a}_B$ is accretive, i.e.\ $\Re\mf{a}_B[u,u] \ge 0$ for all $u\in H^1(\Omega)$. 
    Observe that
	\begin{equation*}
		\Re\mf{a}_B[u,u] \ge \alpha \|\nabla u\|^2_2 + \Re\braket{B\gamma(u),\gamma(u)}_{\partial\Omega}
	\end{equation*}
	for all $u\in \dom(\mf{a}_B) = H^1(\Omega)$, where $\alpha > 0$ is the ellipticity constant from~\ref{setting:main:itm:unif-elliptic}. It follows that $\mf{a}_B$ is accretive if the form $\mf{b}$ is accretive. However, the latter holds if and only if $B+B^*$ is positive semi-definite, due to the identity
	\begin{equation*}
		\braket{(B+B^*)f, f}_{\partial\Omega} 
        = 
        \braket{Bf,f}_{\partial\Omega}+\braket{f,Bf}_{\partial\Omega} 
        = 
        2\Re\braket{Bf,f}_{\partial\Omega} 
        = 
        2\Re \mathfrak{b} [f,f]
	\end{equation*}
	for all $f\in L^2(\partial\Omega)$.

    The accretivity of $\mf{a}_B$ implies
    \begin{equation*}
		\frac{d}{dt}\|e^{-tL_B}u\|^2_2 = -2 \Re\mf{a}_B[e^{-tL_B}u,e^{-tL_B}u] \le 0
	\end{equation*}
	for all $u\in \dom(L_B)$ and all $t\ge 0$. Hence $\|e^{-tL_B}u\|^2_2 \le \|u\|^2_2$ for all $u\in \dom(L_B)$ and all $t\ge 0$. Since $\dom(L_B)$ is dense in $L^2(\Omega)$ (see part~\ref{prop:basic:itm:operator}), we conclude that $(e^{-tL_B})_{t\ge 0}$ is contractive on $L^2(\Omega)$.

    \ref{prop:basic:itm:adjoint} 
    Since $A$ has real coefficients, the adjoint form $\mf{a}_B^*$, defined by $a_B^*[u,v]\coloneqq\overline{\mf{a}_B[v,u]}$ for all $u,v\in\dom(\mf{a}_B)$, is given by
    \begin{equation*}
       \mf{a}^*_B[u,v] = \int_\Omega A^\top \nabla u\cdot\overline{\nabla v}\,dx + \braket{B^*\gamma(u),\gamma(v)}_{\partial\Omega}.
    \end{equation*}
    Using the notation of Setting~\ref{setting:main}\ref{setting:main:itm:A-dependence}, we deduce that $(-L_{A,B})^* = -L_{A^\top, B^*}$. 

    \ref{prop:basic:itm:self-adjoint} This follows immediately from~\ref{prop:basic:itm:adjoint}.

    \ref{prop:basic:itm:real}
    It follows from the characterisation of real semigroups in \cite[Proposition~2.5(1) and~(3)]{Ouh} that $(e^{-tL_B})_{t \ge 0}$ is real if and only if $\mathfrak{a}_B[u,v] \in \RR$ for all real-valued $u,v \in H^1(\Omega)$. 
    We use this to prove the claimed equivalence: 
    
    ``$\Rightarrow$'' 
    If $B$ is real, then clearly $\mathfrak{a}_B[u,v] \in \RR$ for all real-valued $u,v \in H^1(\Omega)$ and thus the semigroup is real. 

    ``$\Leftarrow$'' 
    Assume conversely that the semigroup is real. 
    For all real-valued $u,v \in H^1(\Omega)$, one then has $\mathfrak{a}_B[u,v] \in \RR$ and thus $\braket{B\gamma(u), \gamma(v)}_{\partial \Omega} \in \RR$.
    As the image of the real-valued functions in $H^1(\Omega)$ under $\gamma$ is dense in the real-valued functions in $L^2(\partial \Omega)$ (this follows for instance from \cite[pp.\ 98--99 and Theorem 3.37]{McLean}), we conclude that $\braket{Bf, g}_{\partial \Omega} \in \RR$ for all real-valued $f,g \in L^2(\partial \Omega)$. 
    Thus, $B$ is real.
\end{proof}

In Setting~\ref{setting:main}, for given $f\in L^2(\Omega)$ and $\lambda\in\RR$, we say that $u\in H^1(\Omega)$ is a weak solution of the generalised Robin boundary value problem
\begin{equation}
	\label{eq:bvp-lambda}
	\begin{aligned}
		\lambda u - \dv(A\nabla u) &= f \qquad\text{in }\Omega \\
		\nu\cdot A\nabla u + B\gamma(u) &= 0 \qquad\text{on } \partial\Omega
	\end{aligned}
\end{equation}
if it holds that
\begin{equation*}
    \mf{a}_B[u,v] + \lambda \braket{u,v}_\Omega = \braket{f,v}_\Omega \qquad\forall\, v\in H^1(\Omega).
\end{equation*}
From Proposition~\ref{prop:basic}\ref{prop:basic:itm:H1-elliptic} and the Lax-Milgram theorem, we immediately obtain a well-posedness result.

\begin{corollary}
\label{cor:weak-sol}
    Let the conditions of Setting~\ref{setting:main} be satisfied. 
    Then there exists $\lambda_0>0$ (depending only on $A$, $B$ and $\Omega$)
    such that for every $f\in L^2(\Omega)$ and every $\lambda\ge\lambda_0$, the boundary value problem~\eqref{eq:bvp-lambda}
	has a unique weak solution $u\in H^1(\Omega)$.
\end{corollary}

\begin{remark}
	We briefly remark on higher regularity of solutions. Since the weak solution of~\eqref{eq:bvp-lambda} belongs to $H^1(\Omega)$, we have $B\gamma(u)\in L^2(\partial\Omega)$, and thus the problem of higher regularity reduces to the study of the inhomegeneous Neumann problem
	\begin{equation*}
		\begin{aligned}
			\lambda u - \dv(A\nabla u) &= f \qquad\text{in }\Omega \\
			\nu\cdot A\nabla u &= g \qquad\text{on } \partial\Omega
		\end{aligned}
	\end{equation*}
	with $f\in L^2(\Omega)$ and $g\in L^2(\partial\Omega)$. It is well-known that even in the case of the Laplacian (i.e.\ $A(x) = \id$ for all $x\in\Omega$), $f\in C^\infty(\overline{\Omega})$ and $g\equiv 0$, one cannot expect the `usual' result $u\in H^2(\Omega)$. Counterexamples can be constructed on suitable cones in $\RR^2$ --- some details are given in~\cite[Theorem 1.4.5.3]{Gr}. Thus the regularity problem is highly non-trivial, and has been extensively studied, notably in the works of Jerison and Kenig~\cite{JK}, Fabes et al.~\cite{FMM}, and Savar\'{e}~\cite{Sav}. In the special case of the Laplacian, one has the precise result
    \begin{equation*}
        u\in H^{3/2}_\Delta(\Omega) \coloneqq \{u\in H^{3/2}(\Omega) ~\colon \Delta u\in L^2(\Omega)\};
    \end{equation*}
    see~\cite[Proposition 2.4]{DKP} for a proof.
	
	On the other hand, if we consider $f, g$ with sufficiently high integrability, then Nittka has shown in~\cite[Proposition 3.6]{N11} that the solution to the Neumann problem belongs to the H\"{o}lder space $C^{0,\theta}(\overline{\Omega})$ for some $\theta > 0$.
\end{remark}

We close this section with the following characterisation of positivity of the semigroup $(e^{-tL_B})_{t\ge 0}$, which we will need later. 
In the case of the Laplacian, it can be found in~\cite[Proposition 11.7.1]{GTh}. 
The result generalises to $L_B$ without too much difficulty.

\begin{proposition}
	\label{prop:positivity}
	In Setting~\ref{setting:main}, the following assertions are equivalent:
	\begin{enumerate}[\upshape(i)]
		\item\label{prop:positivity:itm:LB} 
        The semigroup $(e^{-tL_B})_{t\ge 0}$ on $L^2(\Omega)$ is positive.
        
		\item\label{prop:positivity:itm:B} 
        The semigroup $(e^{-tB})_{t\ge 0}$ on $L^2(\partial\Omega)$ is positive.
	\end{enumerate}
\end{proposition}

\begin{proof}
	We use two ingredients in the proof. 
    Firstly, we recall the Beurling--Deny criterion~\cite[Theorem 2.6]{Ouh}. 
    Consider a $\sigma$-finite measure space $(M,\mu)$ and a $C_0$-semigroup $(e^{-tC})_{t\ge 0}$ on $L^2(M,\mu)$ whose generator $-C$ is associated to a sesquilinear form $\mathfrak{c} \colon \dom(\mathfrak{c}) \times \dom(\mathfrak{c}) \to \CC$. 
    The Beurling--Deny criterion states that the semigroup $(e^{-tC})_{t\ge 0}$ is positive if and only if it is real and for every real-valued function $u \in \dom(\mathfrak{c})$ one has $u^+ \in \dom(\mathfrak{c})$ and $\mathfrak{c}[u^+, u^-] \le 0$.

    Secondly, we use that for every real-valued $u \in H^1(\Omega)$, Stampacchia's lemma~\cite[Lemma 7.6]{GT} gives $\int_\Omega A\nabla u^+ \cdot \nabla u^- \,dx = 0$ and therefore
    \begin{equation}
        \label{eq:prop:positivity:eq:proof}
		\mf{b}[\gamma(u)^+,\gamma(u)^-] 
        = 
        \braket{B\gamma(u^+), \gamma(u^-)}_{\partial\Omega} 
        =
        \mf{a}_B[u^+, u^-]
        .
	\end{equation}
    Now we prove the claimed equivalence.
    
    ``\ref{prop:positivity:itm:LB} $\Rightarrow$~\ref{prop:positivity:itm:B}'' 
    Assume that $(e^{-tL_B})_{t\ge 0}$ is positive. 
    Then it is, in particular, real and hence $B$ is real according to Proposition~\ref{prop:basic}\ref{prop:basic:itm:real}. 
    Thus, $(e^{-tB})_{t \ge 0}$ is also real.
    
    Moreover, the Beurling--Deny criterion for $\mathfrak{a}_B$ and the equality~\eqref{eq:prop:positivity:eq:proof} imply that $\mf{b}[\gamma(u)^+,\gamma(u)^-] = \mf{a}_B[u^+, u^-] \le 0$ for every real-valued function $u \in H^1(\Omega)$. 
    Since the image of the real-valued functions in $H^1(\Omega)$ under $\gamma$ is dense in the real-valued functions in $L^2(\partial \Omega)$ (as follows for instance from \cite[pp.\ 98--99 and Theorem 3.37]{McLean}) and $\mf{b} \colon L^2(\partial \Omega) \times L^2(\partial \Omega) \to \CC$ is continuous, it follows that $\mf{b}[f^+,f^-] \le 0$ for all real-valued $f \in L^2(\partial \Omega)$. 
    Thus, the Beurling--Deny criterion applied to the form $\mf{b}$ shows that $(e^{-tB})_{t \ge 0}$ is positive.
    
    ``\ref{prop:positivity:itm:B} $\Rightarrow$~\ref{prop:positivity:itm:LB}''
    Let $(e^{-tB})_{t\ge 0}$ be positive. 
    Then it is, in particular, real, and thus $B$ is real as well. 
    According to Proposition~\ref{prop:basic}\ref{prop:basic:itm:real} this shows that $(e^{-tL_B})_{t \ge 0}$ is also real. 
    Moreover, the Beurling--Deny criterion for $\mf{b}$ gives $\mf{b}[f^+,f^-] \le 0$ for all real-valued $f \in L^2(\partial \Omega)$. 
    Thus, it follows from~\eqref{eq:prop:positivity:eq:proof} that $\mf{a}_B[u^+, u^-] \le 0$ for all real-valued $u \in H^1(\Omega)$. 
    Hence, the Beurling--Deny criterion for $\mf{a}_B$ implies the positivity of $(e^{-tL_B})_{t \ge 0}$.
\end{proof}

\begin{remark}
	\label{rmk:classical-bc}
    Proposition~\ref{prop:positivity} applies in particular to the classical case of local boundary conditions, i.e.\ when $B$ acts as multiplication by a real-valued function $\beta \in L^\infty(\partial\Omega)$. 
    In this case, we write $L_\beta$ instead of $L_B$. 
    Obviously the multiplication semigroup $(e^{-t\beta})_{t\ge 0}$ on $L^2(\partial\Omega)$ is positive and so $(e^{-tL_\beta})_{t\ge 0}$ is positive by Proposition~\ref{prop:positivity}. 

    The positivity of $(e^{-tL_\beta})_{t\ge 0}$ is, of course, known.
    For $\beta \ge 0$ it has been established in~\cite[Theorem 4.9]{AtE97} and~\cite[Proposition 8.1]{Dan00}, 
    and for real-valued $\beta$ is has been proved in~\cite{Dan09}.
\end{remark}

\section{Ultracontractivity}
\label{sec:ultra}

In Setting~\ref{setting:main}, if $B \in \mathcal{L}(L^2(\partial\Omega))$ is given by multiplication with a real-valued function $\beta \in L^\infty(\partial\Omega)$, it is well-known that the semigroup $(e^{-tL_\beta})_{t\ge 0}$ is \emph{ultracontractive}, which means
\begin{equation*}
	e^{-tL_\beta}(L^2(\Omega)) \subseteq L^\infty(\Omega) \qquad \forall \, t>0.
\end{equation*}
Here we are using the notation $L_\beta = L_B$ introduced in Remark~\ref{rmk:classical-bc}.

We now discuss sufficient conditions for the non-local generalisation $-L_B$ to generate an ultracontractive semigroup. 
For this purpose, it is natural to assume that the operator $B$ appearing in the boundary conditions acts boundedly on $L^1(\partial\Omega)$ and $L^\infty(\partial\Omega)$. 
Since this property plays an important role in what follows, let us discuss it in more detail.

Let $(M,\mu)$ be a finite measure space and let $T \in \mathcal{L}\big(L^2(M,\mu)\big)$. 
We say that \emph{$T$ acts boundedly on $L^1(M,\mu)$}, or more briefly that \emph{$T$ acts boundedly on $L^1$}, if $T$ extends to a bounded linear operator on $L^1(M,\mu)$; 
in this case we write $\norm{T}_{L^1 \to L^1}$ for the operator norm of the extension of $T$ to $L^1(M,\mu)$.
Similarly, we say that \emph{$T$ acts boundedly on $L^\infty(M,\mu)$}, or more briefly that \emph{$T$ acts boundedly on $L^\infty$}, if $T$ maps $L^\infty(M,\mu)$ into itself.
By the closed graph theorem, the latter property implies that the restriction of $T$ to $L^\infty(M,\mu)$ is then bounded from $L^\infty(M,\mu)$ to $L^\infty(M,\mu)$, so the terminology \emph{acts boundedly on $L^\infty$} is indeed justified. In this case we write $\norm{T}_{L^\infty \to L^\infty}$ for the operator norm of the restriction of $T$ to $L^\infty(M,\mu)$.

We collect a number of useful results on those properties in the following lemma.

\begin{lemma}
    \label{lem:op-L_1-L_infty}
    Let $(M,\mu)$ be a finite measure space and $T \in \mathcal{L}\big(L^2(M,\mu)\big)$. 
    \begin{enumerate}[\upshape (i)]
        \item \label{lem:op-L_1-L_infty:itm:adjoint} 
        The operator $T$ acts boundedly on $L^1$ if and only if $T^*$ acts boundedly on $L^\infty$. 
        In this case, $\norm{T}_{L^1 \to L^1} = \norm{T^*}_{L^\infty \to L^\infty}$.
        
        \item \label{lem:op-L_1-L_infty:itm:adjoint-2}
        The operator $T$ acts boundedly on $L^\infty$ if and only if $T^*$ acts boundedly on $L^1$.
        In this case, $\norm{T^*}_{L^1 \to L^1} = \norm{T}_{L^\infty \to L^\infty}$.
        
        \item \label{lem:op-L_1-L_infty:itm:modulus}
        Assume that $T$ acts boundedly on $L^1$ and $L^\infty$. 
        Then there exists a positive operator $S \in \mathcal{L}(L^2(M,\mu))$ that also acts boundedly on $L^1$ and $L^\infty$ and satisfies
        \begin{equation}
            \label{lem:modulus-op:eq:estimate}
            \modulus{Tf} \le S\modulus{f} \qquad \forall\, f\in L^2(\Omega,\mu).
        \end{equation}
    \end{enumerate}
\end{lemma}

\begin{proof}
    \ref{lem:op-L_1-L_infty:itm:adjoint} 
    Assume that $T$ acts boundedly on $L^1(M,\mu)$. 
    We first show that $T^* h \in L^\infty(M,\mu)$ for all $h\in L^\infty(M,\mu)$. 
    Indeed, if this is not the case, then there exist $h_0 \in L^\infty(M,\mu)$ and for each $n\in\NN$ a subset $A_n \subset M$ with $\mu(A_n) > 0$ such that $(T^* h_0)(\omega) \ge n$ for almost all $\omega \in A_n$. 
    Define the functions
    \[
        g_n \coloneqq \frac{1}{\mu(A_n)}\one_{A_n}\overline{\operatorname{sgn}(T^* h_0)} 
        \in L^2(M,\mu)
        \qquad\forall\, n\in\NN;
    \]
    clearly $\|g_n\|_{L^1} = 1$ for all $n\in\NN$. 
    It follows that
    \[
        \braket{T^* h_0, g_n} = \frac{1}{\mu(A_n)}\int_{A_n} \modulus{T^* h_0} \,d\mu \ge n \qquad\forall\, n\in\NN.
    \]
    On the other hand,
    \[
        \modulus{\braket{T^* h_0, g_n}} = \modulus{\braket{h_0, Tg_n}} \le \|h_0\|_{L^\infty} \|T\|_{L^1 \to L^1} < \infty,
    \]
    which leads to a contradiction. Thus $T^*(L^\infty(M,\mu)) \subseteq L^\infty(M,\mu)$ as claimed. Moreover
    \[
        \modulus{\braket{T^* h, g}} \le \|h\|_{L^\infty} \|g\|_{L^1} \|T\|_{L^1 \to L^1} \qquad\forall\, h\in L^\infty(M,\mu), g\in L^2(M,\mu).
    \]
    The above estimate then extends by density to all $g\in L^1(M,\mu)$, and therefore $\|T^*\|_{L^\infty \to L^\infty} \le \|T\|_{L^1 \to L^1}$.
    
    Conversely, suppose $T^*$ acts boundedly on $L^\infty(M,\mu)$, and let $f\in L^2(M,\mu)$. Then $Tf \in L^1(M,\mu)$, since $M$ has finite measure. Choose $h\in L^\infty(M,\mu) = (L^1(M,\mu))'$ with $\|h\|_{L^\infty} = 1$ such that $\modulus{\braket{h,Tf}} = \|Tf\|_{L^1}$. Consequently
    \[
        \|Tf\|_{L^1} = \modulus{\braket{T^* h, f}} \le \|T^* h\|_{L^\infty} \|f\|_{L^1} \le \|T^* h\|_{L^\infty \to L^\infty} \|f\|_{L^1}.
    \]
    Since $L^2(M,\mu)$ is dense in $L^1(M,\mu)$, it follows that $T$ extends to a bounded linear operator on $L^1(M,\mu)$ with $\|T\|_{L^1\to L^1} \le \|T^*\|_{L^\infty \to L^\infty}$.

    \ref{lem:op-L_1-L_infty:itm:adjoint-2}
    This follows by applying part~\ref{lem:op-L_1-L_infty:itm:adjoint} to $T^*$ in place of $T$, and by using the fact that $T^{**} = T$.

    \ref{lem:op-L_1-L_infty:itm:modulus} 
    Let $T_1$ denote the extension of $T$ to $L^1(M,\mu)$. Under the given the assumptions on $T$, a classical result~\cite[Lemma VIII.6.4]{DunfordSchwartz-I} shows that there exists a positive operator $S_1$ on $L^1(M,\mu)$ such that
    \[
        \modulus{T_1 g} \le S_1 \modulus{g} \qquad\forall\, g\in L^1(M,\mu),
    \]
    and such that $S_1(L^\infty(M,\mu)) \subseteq L^\infty(M,\mu)$ and $\norm{S_1}_{L^\infty \to L^\infty} < \infty$. 
    By the Riesz--Thorin interpolation theorem, the restriction $S \coloneqq S_1 \vert_{L^2(M,\mu)}$ is a bounded operator on $L^2(M,\mu)$, and clearly $S$ has all the required properties.
\end{proof}

As we will see, the assumption that the operator $B$ acts boundedly on $L^1(\partial\Omega)$ and $L^\infty(\partial\Omega)$ ensures that each of the semigroup operators $e^{-tL_B}$ acts boundedly on $L^1(\Omega)$ and $L^\infty(\Omega)$, which then enables us to apply existing results on ultracontractivity.
The following theorem is the main result of this section.

\begin{theorem}[Ultracontractivity]
    \label{thm:ultra}
    In Setting~\ref{setting:main}, assume that the operator $B\in\mathcal{L}(L^2(\partial \Omega))$ acts boundedly on $L^1$ and $L^\infty$.
    Then the semigroup $(e^{-tL_B})_{t\ge 0}$ satisfies $e^{-tL_B} (L^2(\Omega)) \subseteq L^\infty(\Omega)$ for all $t>0$, and there exist constants $c,\mu>0$ such that
    \begin{equation}
	   \label{eq:ultra-estimate}
	   \|e^{-tL_B}\|_{L^2 \to L^\infty} \le ct^{-\mu/4} \qquad \forall \, t \in (0,1]
    \end{equation}
    and such that the adjoint operators $e^{-tL_B^*}=(e^{-tL_B})^*$ also satisfy~\eqref{eq:ultra-estimate}.
\end{theorem}

Before giving the proof, we mention two simple examples of operators $B\in\mathcal{L}(L^2(\partial \Omega))$ that satisfy the assumptions of Theorem~\ref{thm:ultra}. Further examples will be exhibited later in Examples~\ref{exa:operators-B-ev-pos}.

\begin{examples}
    \label{exa:boundary-op}
    \begin{enumerate}[(a)]
        \item 
        Let the operator $B\in\mathcal{L}(L^2(\partial\Omega))$ be given by multiplication with a function $\beta \in L^\infty(\partial\Omega)$. Clearly, such an operator acts boundedly on $L^1$ and $L^\infty$. In this setting, we write $L_\beta$ instead  $L_B$.
        When $\beta$ is real-valued, we recover the classical Robin boundary conditions mentioned already in Remark~\ref{rmk:classical-bc}.

        \item 
        Another simple class of operators $B\in\mathcal{L}(L^2(\partial\Omega))$ that satisfies the assumptions of Theorem~\ref{thm:ultra} are kernel operators
        \[
            Bf \coloneqq \int_{\partial\Omega} k(\,\cdot\, ,y)f(y)\,d\sigma(y) \qquad\forall\,f\in L^2(\partial\Omega),
        \]
        where $k\in L^\infty(\partial\Omega \times \partial\Omega)$ and $\sigma$ is the $(d-1)$-dimensional Hausdorff measure on $\partial\Omega$.
        Since $\Omega$ is a bounded Lipschitz domain, we have $\sigma(\partial\Omega) < \infty$, thus $B$ acts boundedly on $L^1$ and $L^\infty$.
    \end{enumerate}
\end{examples}

The proof of Theorem~\ref{thm:ultra} requires some preparation. The most non-trivial part is to show that the semigroup operators $e^{-tL_B}$ extend to bounded operators on $L^1$ and $L^\infty$ with uniform bounds for $t\in [0,1]$. Since the semigroup $(e^{-tL_B})_{t\ge 0}$ need not be positive, this will be achieved by constructing a dominating positive semigroup, for which it is easier to obtain the desired uniform bounds.

We begin with a simple but general criterion for $L^\infty$ boundedness for positive semigroups.

\begin{proposition}
    \label{prop:Linfty-bounded-pos-sgrp}
    Let $(M,\mu)$ be a finite measure space and let $(e^{-tA})_{t\ge 0}$ be a positive $C_0$-semigroup on $L^2(M,\mu)$. 
    Suppose there exists a real-valued function $u\in L^\infty(M,\mu)\cap \dom(A)$ satisfying the following properties:
    \begin{enumerate}[\upshape (i)]
        \item\label{prop:Linfty-bounded-pos-sgrp:it:strong-pos} 
        $u \ge \delta\one$ for some constant $\delta>0$;
        
        \item\label{prop:Linfty-bounded-pos-sgrp:it:A-estimate}
        There exists $\nu\in\RR$ in the resolvent set of $-A$ such that $(\nu + A)u \ge 0$.
    \end{enumerate}
    Then $e^{-tA}$ restricts to a bounded operator on $L^\infty(M,\mu)$ for every $t\ge 0$, and one has $\sup_{t \in [0,1]} \|e^{-tA}\|_{L^\infty\to L^\infty} < \infty$.
\end{proposition}

\begin{proof}
    By increasing $\nu$ if necessary we may assume that $\nu > \omega_0(-A)$; 
    note that increasing $\nu$ does not affect the assumption~\ref{prop:Linfty-bounded-pos-sgrp:it:A-estimate} since $u$ is positive due to~\ref{prop:Linfty-bounded-pos-sgrp:it:strong-pos}. 
    Set $w \coloneqq (\nu + A)u \ge 0$. 
    The Laplace transform representation of the resolvent yields $u = R(\nu,-A)w = \int_0^\infty e^{-\nu s}e^{-sA}w \,ds$. Consequently
    \[
        e^{-tA}u = \int_0^\infty e^{-\nu s} e^{-(t+s)A}w\, ds = e^{\nu t}\int_t^\infty e^{-\nu\tau} e^{-\tau A}w \, d\tau \le e^{\nu t} u
    \]
    for all $t\ge 0$, where the inequality is due to the positivity of the semigroup and of $w$. 
    For every $f\in L^\infty(M,\mu)$ with $\|f\|_\infty \le 1$,
    assumption~\ref{prop:Linfty-bounded-pos-sgrp:it:strong-pos} implies $\modulus{f}\le \delta^{-1} u$ and hence
    \[
        \modulus{e^{-tA}f} \le e^{-tA}\modulus{f} \le \delta^{-1} e^{\nu t}u
    \]
    for all $t\ge 0$. Since $u\in L^\infty(\Omega)$, the assertions of the proposition follow readily.
\end{proof}

Next, we require a regularity result for the case of local boundary conditions.

\begin{lemma}
    \label{lem:robin-local}
    In Setting~\ref{setting:main}, assume that the operator $B\in\mathcal{L}(L^2(\Omega))$ is given by multiplication with a real-valued function $\beta\in L^\infty(\partial\Omega)$;
    we thus write $L_\beta \coloneqq L_B$.

    For all sufficiently large $\lambda>0$, the function $u \coloneqq \lambda R(\lambda,-L_\beta)\one$ belongs to $C(\overline{\Omega})$ and satisfies $\frac{1}{2} \le u(x) \le 2$ for all $x\in\overline{\Omega}$.
\end{lemma}

\begin{proof}
    It was shown by Nittka~\cite[Theorem 4.3]{N11} that the part $-L_{\beta,C}$ of $-L_\beta$ in $C(\overline{\Omega})$ generates a positive $C_0$-semigroup on $C(\overline{\Omega})$.
    It follows that $\lambda R(\lambda,-L_\beta)\one = \lambda R(\lambda,-L_{\beta,C})\one$ is a real-valued function in $C(\overline{\Omega})$ and converges to $\one$ with respect to the supremum norm as $\lambda \to \infty$.
    Hence for sufficiently large $\lambda$ we indeed have $\frac{1}{2} \le u(x) \le 2$ for all $x\in\overline{\Omega}$. 
    We emphasise that the inequality makes sense on $\overline{\Omega}$ rather than only on $\Omega$, since $u \in C(\overline{\Omega})$.
\end{proof}

The final key ingredient for the proof of Theorem~\ref{thm:ultra} is the following consequence of Ouhabaz' criterion for domination between semigroups~\cite[Theorem 2.21]{Ouh}.

\begin{proposition}
    \label{prop:domination-robin}
    Let the assumptions of Setting~\ref{setting:main} be satisfied and let $\tilde B \in \mathcal{L}(L^2(\partial\Omega))$. 
    Consider the following assertions:
    \begin{enumerate}[\upshape(i)]
        \item\label{prop:domination-robin:itm:B}
        $\braket{\tilde B\modulus{f},\modulus{g}}_{\partial \Omega} \le \Re \braket{Bf,g}_{\partial \Omega} $ for all $f,g \in L^2(\partial \Omega)$.
        
        \item\label{prop:domination-robin:itm:domination} 
        $\modulus{e^{-tL_B}f} \le e^{-tL_{\tilde B}}\modulus{f} $ for all $f\in L^2(\Omega)$ and all $t \ge 0$.
        
        \item\label{prop:domination-robin:itm:forms} 
        $\tilde B$ is real and $\mf{a}_{\tilde B}[\modulus{v},\modulus{w}] \le \Re\mf{a}_B[v,w]$ for all $v,w\in H^1(\Omega)$ such that $v\overline{w}\ge 0$.
    \end{enumerate}
    Then \ref{prop:domination-robin:itm:B}~$\Rightarrow$~\ref{prop:domination-robin:itm:domination} and \ref{prop:domination-robin:itm:domination}~$\Leftrightarrow$~\ref{prop:domination-robin:itm:forms}. 
    If one of those three assertions holds, then $(e^{-t\tilde B})_{t \ge 0}$ and $(e^{-tL_{\tilde B}})_{t \ge 0}$ are positive.
\end{proposition}

\begin{proof}
    Let us first show that each of the three conditions~\ref{prop:domination-robin:itm:B}--\ref{prop:domination-robin:itm:forms} implies the positivity of $(e^{-tL_{\tilde B}})_{t \ge 0}$ (equivalently, the positivity of $(e^{-t\tilde B})_{t \ge 0}$, see Proposition~\ref{prop:positivity}). 

    First, let~\ref{prop:domination-robin:itm:B} hold. 
    If $f,g \in L^2(\partial \Omega)$ satisfy $f,g \ge 0$, then~\ref{prop:domination-robin:itm:B} can be applied to $f$ and $g$, and also to $-f$ and $g$, to obtain $\braket{\tilde B f, g}_{\partial \Omega} \le 0$. 
    Hence, $-\tilde B$ is a positive operator and thus $(e^{-t\tilde B})_{t \ge 0}$ is positive. 

    Next, observe that~\ref{prop:domination-robin:itm:domination} immediately implies the positivity of $(e^{-tL_{\tilde B}})_{t \ge 0}$. 

    Finally, let~\ref{prop:domination-robin:itm:forms} hold. 
    Since $\tilde B$ is real, the semigroup $(e^{-tL_{\tilde B}})_{t \ge 0}$ is real by Proposition~\ref{prop:basic}\ref{prop:basic:itm:real}. 
    Now let $v \in H^1(\Omega)$ be real-valued. 
    Note that $v^+ v^- = 0$. 
    By applying~\ref{prop:domination-robin:itm:forms} to $v^+$ and $v^-$, and then to $-v^+$ and $v^-$, we obtain $\mf{a}_{\tilde B}[v^+,v^-] \le 0$. 
    Thus, the Beurling--Deny criterion \cite[Theorem~2.6(1) and~(4)]{Ouh} implies that $(e^{-tL_{\tilde B}})_{t \ge 0}$ is positive.
    \smallskip 
    
    Now we show the claimed implications.
    
    ``\ref{prop:domination-robin:itm:domination} $\Leftrightarrow$ \ref{prop:domination-robin:itm:forms}'' 
    As both conditions imply the positivity of $(e^{-tL_{\tilde B}})_{t \ge 0}$, their equivalence follows from the characterisation of domination between semigroups given in \cite[Theorem~2.21(1) and~(3)]{Ouh}. 
    Note that the theorem in this reference is only stated for accretive forms, but by rescaling the semigroups one immediately gets the same result for forms that satisfy the ellipticity estimate in Proposition~\ref{prop:basic}\ref{prop:basic:itm:H1-elliptic}.
    
    ``\ref{prop:domination-robin:itm:B} $\Rightarrow$ \ref{prop:domination-robin:itm:domination}'' 
    Assume~\ref{prop:domination-robin:itm:B}. 
    We actually show that this implies~\ref{prop:domination-robin:itm:forms}.  

    At the beginning of the proof we have shown that~\ref{prop:domination-robin:itm:B} implies $-\tilde B \ge 0$.
    In particular, $\tilde B$ is real.
    
    Now let $v,w \in H^1(\Omega)$ satisfy $v \overline{w} \ge 0$.
    The semigroup $(e^{-tL_0})_{t \ge 0}$ associated to the form $\mf{a}_0$ is positive (Proposition~\ref{prop:positivity}) and thus satisfies $\modulus{e^{-tL_0}h} \le e^{-tL_0}\modulus{h}$ for all $t \ge 0$ and all $h \in L^2(\Omega)$. 
    Thus, if we apply the equivalence of~\ref{prop:domination-robin:itm:domination} and~\ref{prop:domination-robin:itm:forms} (to the special case $B = \tilde B = 0$), we obtain $\mf{a}_0[\modulus{v},\modulus{w}] \le \Re \mf{a}_0[v,w]$. 
    Since $\gamma(\modulus{v}) = \modulus{\gamma(v)}$, and likewise for $w$,~\ref{prop:domination-robin:itm:B} gives $\braket{\tilde B\gamma(\modulus{v}),\gamma(\modulus{w})}_{\partial \Omega} \le \Re \braket{B\gamma(v),\gamma(w)}_{\partial \Omega}$ and hence,
    \begin{align*}
        \mf{a}_{\tilde B}[\modulus{v}, \modulus{w}] 
        & =
        \mf{a}_{0}[\modulus{v}, \modulus{w}] 
        + 
        \braket{\tilde B\gamma(\modulus{v}),\gamma(\modulus{w})}_{\partial \Omega}
        \\ 
        & \le 
        \Re \mf{a}_0[v,w]  +  \Re \braket{B\gamma(v),\gamma(w)}_{\partial \Omega} 
        = 
        \Re \mf{a}_B[v,w]. 
        \qedhere
    \end{align*}
\end{proof}

In the preceding proposition, positivity of the semigroup $(e^{-tL_{\tilde B}})_{t \ge 0}$ is a consequence of each of the assertions~\ref{prop:domination-robin:itm:B}--\ref{prop:domination-robin:itm:forms}, while it is an assumption in~\cite[Theorem~2.21]{Ouh}. 
This difference is due to the fact that both forms $\mf{a}_B$ and $\mf{a}_{\tilde B}$ have the same domain in Proposition~\ref{prop:domination-robin}. 

All the ingredients are now in place to prove our ultracontractivity result stated in Theorem~\ref{thm:ultra}.

\begin{proof}[Proof of Theorem~\ref{thm:ultra}]
    We proceed in several steps.
    \begin{description}[font=\normalfont\itshape, labelindent=0cm, leftmargin=0.4cm]
        \item[Step~1]
        According to Lemma~\ref{lem:op-L_1-L_infty}\ref{lem:op-L_1-L_infty:itm:modulus}, there is  a positive operator $-B_1 \in \mathcal{L}(L^2(\partial \Omega))$ that acts boundedly on $L^1$ and $L^\infty$ and satisfies 
        \begin{align*}
            \modulus{Bf} \le -B_1 \modulus{f} 
            \qquad 
            \forall\, f \in L^2(\partial \Omega)
            .
        \end{align*}

        \item[Step~2]
        Now we use properties of local Robin boundary conditions as an auxiliary tool:
        we construct functions $0 \le -\beta \in L^\infty(\partial \Omega)$, $0 \le u \in \dom(L_{\beta})$, and a number $\lambda > 0$ with the following properties:
        \begin{enumerate}[(1), leftmargin=1cm]
            \item \label{thm:ultra:step1:itm:bnds}
            $u \in C(\overline{\Omega})$ and $\frac{1}{2} \le u(x) \le 2$ for all $x \in \overline{\Omega}$;

            \item \label{thm:ultra:step1:itm:pos}
            $\lambda u + L_\beta u \ge 0$;
            
            \item \label{thm:ultra:step1:itm:boundary}
            $-B_1 \gamma(u) \le -\beta \gamma(u)$.
        \end{enumerate}
        To this end, choose $\beta \coloneqq -4 \norm{B_1}_{L^\infty \to L^\infty} \one_{\partial \Omega}$, which is possible since $B_1$ acts boundedly on $L^\infty(\partial \Omega)$ by Step~1. 
        Then Lemma~\ref{lem:robin-local} shows that there exists a $\lambda > 0$ in the resolvent set of $-L_\beta$ such that $u \coloneqq \lambda R(\lambda, -L_\beta)\one$ satisfies property~\ref{thm:ultra:step1:itm:bnds}. 
        Property~\ref{thm:ultra:step1:itm:pos} holds, since we have 
        \begin{align*}
            \lambda u + L_\beta u = \lambda \one \ge 0
        \end{align*}
        by the definition of $u$ and the positivity of $\lambda$.
        Property~\ref{thm:ultra:step1:itm:boundary} now follows from~\ref{thm:ultra:step1:itm:bnds} and the definition of $\beta$, since 
        \begin{align*}
            -B_1 \gamma(u) 
            \le 
            -2B_1 \one_{\partial \Omega} 
            \le 
            2 \norm{B_1}_{L^\infty \to L^\infty} \one_{\partial \Omega} 
            \le 
            4 \norm{B_1}_{L^\infty \to L^\infty} \gamma(u) 
            = 
            -\beta \gamma(u)
            ,
        \end{align*}
        where we used the positivity of $-B_1$ for the first inequality.

        \item[Step~3]
        Next, we construct an operator $B_2 \in \mathcal{L}(L^2(\partial \Omega))$ that satisfies
        \begin{enumerate}[(1), start=4, leftmargin=1cm]
            \item \label{thm:ultra:step1:itm:boundary-equ} 
            $-B_2 \gamma(u) = -\beta \gamma(u)$;

            \item \label{thm:ultra:step1:itm:dom-B}
            $0 \le -B_1 \le -B_2$;

            \item \label{thm:ultra:step1:itm:change-op}
            $u \in \dom (L_{B_2})$ and $L_{B_2}u = L_\beta u$.
        \end{enumerate}
        Note that~\ref{thm:ultra:step1:itm:boundary-equ} means that $B_2$ improves the inequality in property~\ref{thm:ultra:step1:itm:boundary} to an equality.
        
        To construct $B_2$, consider the function $g \coloneqq -\beta \gamma(u) + B_1 \gamma(u) \in L^2(\partial \Omega)$ and choose $0 \le h \in L^2(\partial \Omega)$ such that $\braket{\gamma(u), h}_{\partial \Omega} = 1$. 
        Note that $g \ge 0$ by property~\ref{thm:ultra:step1:itm:boundary}.
        Then the operator $B_2$ defined by 
        \begin{align*}
            B_2 f 
            \coloneqq 
            B_1 f - \braket{f, h}_{\partial \Omega} g 
            \qquad 
            \forall\, f \in L^2(\partial \Omega)
        \end{align*}
        satisfies both~\ref{thm:ultra:step1:itm:boundary-equ} and~\ref{thm:ultra:step1:itm:dom-B}.
        To show property~\ref{thm:ultra:step1:itm:change-op}, observe that~\ref{thm:ultra:step1:itm:boundary-equ} implies $\mf{a}_{B_2}[u,v] = \mf{a}_{\beta}[u,v] = \braket{L_\beta u, v}_{\Omega}$ for all $v\in H^1(\Omega)$. Hence, indeed $u\in\dom(L_{B_2})$ and $L_{B_2}u = L_{\beta}u$.

        \item[Step~4] 
        We show that the semigroup $(e^{-tL_{B_2}})_{t \ge 0}$ is positive and that 
        \begin{align}
            \label{eq:thm:ultra:eq:dom}
            \modulus{e^{-tL_B}h} \le e^{-tL_{B_2}} \modulus{h} 
            \qquad \forall\, h \in L^2(\Omega).
        \end{align}
        Indeed, the positivity of $(e^{-tL_{B_2}})_{t \ge 0}$ follows from Proposition~\ref{prop:positivity} since $-B_2 \ge 0$. 
        It thus suffices to check condition~\ref{prop:domination-robin:itm:B} in Proposition~\ref{prop:domination-robin} for the forms $\mf{a}_{B_2}$ and $\mf{a}_{B}$ to establish the domination property~\eqref{eq:thm:ultra:eq:dom}. 
        So let $f,g \in L^2(\partial \Omega)$ and observe that 
        \begin{align*}
             - \Re \braket{Bf,g}_{\partial \Omega} 
             \le 
             \modulus{\braket{Bf,g}_{\partial \Omega}} 
             &\le 
             \braket{\modulus{Bf},\modulus{g}}_{\partial \Omega} \\
             &\le 
             \braket{-B_1\modulus{f},\modulus{g}}_{\partial \Omega} 
             \le 
             \braket{-B_2\modulus{f},\modulus{g}}_{\partial \Omega}
             ,
        \end{align*}
        where the third inequality follows from the choice of $B_1$ in Step~1 and the last inequality uses property~\ref{thm:ultra:step1:itm:dom-B}.
        
        \item[Step~5] 
        We show that each operator of the semigroup $(e^{-tL_{B}})_{t \ge 0}$ acts boundedly on $L^\infty(\Omega)$ and that
        \begin{align*}
            \sup_{t \in [0,1]} \norm{e^{-tL_{B}}}_{L^\infty \to L^\infty} < \infty.
        \end{align*}
        To this end, first note that the function $u$ from Step~2 satisfies $\lambda u + L_{B_2} u \ge 0$ as a consequence of properties~\ref{thm:ultra:step1:itm:pos} and~\ref{thm:ultra:step1:itm:change-op}. 
        Thus, assumption~\ref{prop:Linfty-bounded-pos-sgrp:it:A-estimate} from Proposition~\ref{prop:Linfty-bounded-pos-sgrp} is satisfied. 
        Assumption~\ref{prop:Linfty-bounded-pos-sgrp:it:strong-pos} from the same proposition also holds due to the lower estimate for $u$ in property~\ref{thm:ultra:step1:itm:bnds}. 
        Thus, Proposition~\ref{prop:Linfty-bounded-pos-sgrp} shows that $e^{-tL_{B_2}}$ leaves $L^\infty(\Omega)$ invariant for each $t \ge 0$ and that $\sup_{t \in [0,1]} \norm{e^{-tL_{B_2}}}_{L^\infty \to L^\infty} < \infty$. 
        By the domination property~\eqref{eq:thm:ultra:eq:dom}, the same holds for the semigroup generated by $-L_B$.

        \item[Step~6] 
        We show that each operator of the semigroup $(e^{-tL_{B}})_{t \ge 0}$ acts boundedly on $L^1(\Omega)$ and that
        \begin{align*}
            \sup_{t \in [0,1]} \norm{e^{-tL_{B}}}_{L^1 \to L^1} < \infty.
        \end{align*}
        Indeed, until now we have only used the assumption that $B$ acts boundedly on $L^\infty$. 
        Since we also assumed $B$ to act boundedly on $L^1$, its adjoint $B^*$ acts boundedly on $L^\infty$ as well (Lemma~\ref{lem:op-L_1-L_infty}\ref{lem:op-L_1-L_infty:itm:adjoint}).
        Hence, Steps~1--5 also apply to the operator $L_{A^\top,B^*} = (L_{A,B})^* = (L_B)^*$, where we used the notation $L_{A,B} = L_B$ from Setting~\ref{setting:main}\ref{setting:main:itm:A-dependence} and where the first equality follows from Proposition~\ref{prop:basic}\ref{prop:basic:itm:adjoint}.

        Thus, $(e^{-L_B})^* = e^{-tL_B^*}$ acts boundedly on $L^\infty$ for every $t \ge 0$ and one has $\sup_{t \in [0,1]} \norm{(e^{-tL_B})^*}_{L^\infty \to L^\infty} < \infty$. 
        According to Lemma~\ref{lem:op-L_1-L_infty}\ref{lem:op-L_1-L_infty:itm:adjoint-2} this implies that $e^{-tL_B}$ acts boundedly on $L^1$ for each $t \ge 0$ and that 
        \begin{align*}
            \sup_{t \in [0,1]} \norm{e^{-tL_{B}}}_{L^1 \to L^1} < \infty.
        \end{align*}

        \item[Step~7]
        Finally, we show the claimed ultracontractivity. 

        The analyticity of the semigroup implies that $e^{-tL_B}(L^2(\Omega)) \subseteq \dom(-L_B)$ for all $t>0$ and that there is a constant $c_1 > 0$ such that $\|L_B e^{-tL_B}\|_{L^2\to L^2} \le c_1 t^{-1}$ for all $t\in (0,1]$; see for instance~\cite[Proposition 2.1.1]{Lunardi95} for this standard fact. Combining this with Proposition~\ref{prop:basic}\ref{prop:basic:itm:H1-elliptic}, we see that there are constants $c_0,\omega>0$ such that for all $t\in (0,1]$ and $f\in L^2(\Omega)$ with $\|f\|_2 \le 1$,
        \begin{align*}
            c_0\|\nabla e^{-tL_B}f\|^2_{2} &\le \mf{a}_B[e^{-tL_B}f, e^{-tL_B}f] + \omega\|e^{-tL_B}\|^2_2 \\
            &\le \braket{L_B e^{-tL_B}f, e^{-tL_B}f}_\Omega + M^2 \omega \\
            &\le M c_1 t^{-1} + M^2 \omega
            ,
        \end{align*}
        where $M \coloneqq \sup_{t\in [0,1]}\|e^{-tL_B}\|_{L^2\to L^2} < \infty$.
        It is then easy to deduce that
        \begin{equation}
        \label{thm:ultra:eq:L2-to-H1}
            \|e^{-tL_B}\|_{L^2 \to H^1} \le c t^{-1/2} \qquad\forall\, t\in (0,1]
        \end{equation}
        for a suitable constant $c>0$.
        
        In the case of dimension $d=1$, we may use the estimate~\eqref{thm:ultra:eq:L2-to-H1} combined with the continuous embedding $H^1(\Omega) \hookrightarrow L^\infty(\Omega)$ (cf.~\cite[Theorem 8.8]{Bre}) to obtain the desired ultracontractivity estimate~\eqref{eq:ultra-estimate} with $\mu = 2$.

        To obtain ultracontractivity in dimensions $d\ge 2$, we make use of the uniform bounds in Steps~5 and 6. 
        Firstly, we note that the extension of the semigroup $(e^{-tL_B})_{t \ge 0}$ to $L^1(\Omega)$ is also a $C_0$-semigroup. 
        Indeed, for every $f \in L^2(\Omega)$, the orbit map $t \mapsto e^{-tL_B}f$ is continuous from $[0,\infty)$ into $L^2(\Omega)$, and thus also into $L^1(\Omega)$ by the continuous embedding of $L^2(\Omega)$ into $L^1(\Omega)$. 
        The uniform $L^1$ bound from~Step~6 thus gives the $C_0$-property on $L^1(\Omega)$ \cite[Proposition I.5.3]{EN00}. 
        Consequently, $(e^{-tL_B})_{t\ge 0}$ extrapolates to a consistent family of semigroups acting on the $L^p$ scale $1\le p \le \infty$, and is strongly continuous on $L^p$ for all $p\ne\infty$ (cf.~\cite[p.\ 61]{A04}).
        This extrapolation property together with the analyticity of the semigroup on $L^2(\Omega)$ allow us to apply the characterisation~\cite[Theorem 7.3.2]{A04}, in particular the implication (v) $\Rightarrow$ (ii): if the form domain $V \coloneqq \dom(\mf{a}_B)$ satisfies the continuous embedding $V \hookrightarrow L^{2\mu/(\mu-2)}(\Omega)$ for some $\mu>2$, and $V\cap L^1(\Omega)$ is dense in $L^1(\Omega)$, then the semigroup $(e^{-tL_B})_{t\ge 0}$ is ultracontractive with the estimate
	    \begin{equation*}
		      \|e^{-tL_B}\|_{L^2 \to L^\infty} \le ct^{-\mu/4} \quad\,\forall\,t\in (0,1].
	    \end{equation*}
	    In our situation, we have $V=H^1(\Omega)\supseteq C^\infty_c(\Omega)$, and thus the density condition is clearly satisfied. We conclude using the Sobolev embedding theorems~\cite[Theorem 4.12, Part I]{AF}. If $d\ge 3$, we may choose $\mu=d$ since the embedding $V\hookrightarrow L^{2d/(d-2)}$ is valid in bounded Lipschitz domains. Finally, in dimension $d=2$, we have $V\hookrightarrow L^q(\Omega)$ for any $2\le q<\infty$, so we may choose any $2<\mu<\infty$. In each case, it follows that $(e^{-tL_B})_{t\ge 0}$ is ultracontractive with the desired estimate~\eqref{eq:ultra-estimate}, and the proof is complete.
    \qedhere
    \end{description}
\end{proof}


It is natural to ask if one can obtain Gaussian estimates for the semigroup $(e^{-tL_B})_{t\ge 0}$, which would yield ultracontractivity as an immediate consequence. 
This was already investigated in~\cite{GeMiNi14, GeMiNiOu}, where the authors even allow for certain classes of unbounded boundary operators $B$. However, the assumptions in these articles lead to the domination property
\begin{equation}
    \label{eq:neumann-domination}
    \modulus{e^{-tL_B}f} \le e^{t\Delta_N}\modulus{f} \qquad \forall\, t\ge 0, f\in L^2(\Omega)
\end{equation}
where $\Delta_N$ denotes the Neumann Laplacian (in the notation of Setting~\ref{setting:main}, this is the operator $L_0$ associated to the form $\mf{a}_0$ with $B=0$), see in particular~\cite[Theorem 4.4]{GeMiNi14} and~\cite[Theorem 3.3]{GeMiNiOu}. 
This is an extremely strong property, as the following proposition shows.

Let $\emptyset \neq \Omega \subseteq \RR^d$ be open. 
We call a sesquilinear form $\mathfrak{c} \colon H^1(\Omega) \times H^1(\Omega) \to \CC$ \emph{local} if it satisfies $\mathfrak{c}[v,w] = 0$ for all $v,w \in H^1(\Omega)$ with $vw = 0$.

\begin{proposition}
    \label{prop:local-by-dom}
    Assume that the conditions of Setting~\ref{setting:main} are satisfied. 
    \begin{enumerate}[\upshape(i)]
        \item\label{prop:local-by-dom:itm:local} 
        The form $\mf{a}_B$ is local if and only if $B$ is a multiplication operator, meaning that $B$ acts as the multiplication with a function $\beta \in L^\infty(\partial \Omega)$.
        
        \item\label{prop:local-by-dom:itm:dom} 
        Now let $\tilde B \in \mathcal{L}(L^2(\partial\Omega))$ as well, and suppose that 
        $\modulus{e^{-tL_B}f} \le e^{-tL_{\tilde B}}\modulus{f} $ for all $f\in L^2(\Omega)$ and all $t \ge 0$. 
        If $\mf{a}_{\tilde B}$ is local, then so is $\mf{a}_B$.
    \end{enumerate}
\end{proposition}

\begin{proof}
    \ref{prop:local-by-dom:itm:local}
    First observe that the form $\mf{a}_0$ is local as a consequence of Stampacchia's lemma~\cite[Lemma 7.6]{GT}. 
    Thus, if $B$ is a multiplication operator, it immediately follows that $\mf{a}_B$ is local. 

    Conversely, assume now that $\mf{a}_B$ is local.  
    Since $\mf{a}_0$ is local, it follows that $\braket{B\gamma(u)^+, \gamma(u)^-}_{\partial \Omega} = \braket{B\gamma(u^+), \gamma(u^-)}_{\partial \Omega} = 0$ 
    for every real-valued $u \in H^1(\Omega)$. 
    Since the image of all real-valued functions in $H^1(\Omega)$ under the trace operator $\gamma$ 
    is dense in the real-valued functions in $L^2(\partial \Omega)$ (see for instance~\cite[pp.\ 98--99 and Theorem 3.37]{McLean}), 
    it follows that $\braket{Bf^+, f^-}_{\partial \Omega} = 0$ for all real-valued $f \in L^2(\partial \Omega)$. 
    This in turn implies that $\braket{Bf, g}_{\partial \Omega} = 0$ for all $0 \le f,g \in L^2(\partial \Omega)$ with $fg = 0$. 

    If $f,g \in L^2(\partial \Omega)$ are general real-valued functions with $fg = 0$, then $(f^+ + f^-)(g^+ + g^-) = \modulus{fg} = 0)$. 
    Hence, all the products $f^+g^+$, $f^+ g^-$, $f^- g^+$, $f^-g^-$ are zero, so the sesquilinearity of $\mf{a}_B$ yields $\braket{Bf, g}_{\partial \Omega} = 0$. 
    Finally, if $f,g \in L^2(\partial \Omega)$ are complex-valued functions with $fg = 0$, then 
    \begin{align*}
        0 
        = 
        4 
        \modulus{fg} 
        \ge 
        \big(\modulus{\Re f} + \modulus{\Im f}\big) \big(\modulus{\Re g} + \modulus{\Im g}\big) 
        ,
    \end{align*}
    so all the products $\Re f \Re g$, $\Re f \Im g$, $\Im f \Re g$, $\Im f \Im g$ are zero. 
    Thus, the sesqulinearity of $\mf{a}_B$ implies again that $\braket{Bf, g}_{\partial \Omega} = 0$.

    Next, we show that if $f,g \in L^2(\partial \Omega)$ satisfy $fg = 0$, then even $(Bf)g = 0$. 
    To this end, let $s \in L^\infty(\partial \Omega)$ be given by $s(x) = \overline{\operatorname{sign}(Bf(x))}$ for $x \in \partial \Omega$. 
    Since $fs\modulus{g} = 0$, it follows that 
    \begin{align*}
        0 
        = 
        \braket{B f, s \modulus{g}} 
        = 
        \int_{\partial \Omega} \modulus{Bf} \modulus{g} \, d\sigma
        .
    \end{align*}
    Hence, $\modulus{Bf}\modulus{g} = 0$ and so $(Bf)g = 0$, as claimed. 

    Now we can show that $B$ is indeed a multiplication operator. 
    This actually follows by abstract Banach lattice theory (see for instance \cite[Example~2.67]{AB}), 
    but we prefer to include an elementary proof.
    To this end, set $\beta \coloneqq B\one_{\partial \Omega} \in L^2(\partial \Omega)$. 
    If $A \subseteq \partial \Omega$ is measurable, one has $(B\one_A) \one_{A^{\mathrm{c}}} = 0$ and $(B\one_{A^{\mathrm c}}) \one_A = 0$. 
    Hence, 
    \begin{align*}
        \beta \one_A
        = 
        (B \one_A + B\one_{A^{\mathrm c}}) \one_A 
        = 
        B\one_A
        .
    \end{align*} 
    We show that this implies $\beta \in L^\infty(\partial \Omega)$. 
    To this end set $A \coloneqq \{x \in \partial \Omega : \, \modulus{\beta(x)} \ge \norm{B}_{L^2 \to L^2} + 1\}$. 
    We show that the $(d-1)$-dimensional Hausdorff measure $\sigma(A)$ is zero. 
    Indeed, 
    \begin{align*}
        \norm{B}_{L^2 \to L^2} \sigma(A)^{1/2} 
        \ge 
        \norm{B \one_A}_{2,\partial\Omega} 
        = 
        \norm{\beta \one_A}_{2,\partial\Omega} 
        \ge 
        (\norm{B}_{L^2 \to L^2} + 1) \sigma(A)^{1/2}
        ,
    \end{align*}
    which shows that $\sigma(A) = 0$, as claimed. 
    So indeed $\beta \in L^\infty(\partial \Omega)$. 
    
    By a density argument, the fact that $B \one_A = \beta \one_A$ for all measurable $A \subseteq \partial\Omega$ now implies that $Bf = \beta f$ for all $f \in L^2(\partial \Omega)$.

    \ref{prop:local-by-dom:itm:dom}
    Let $v,w \in H^1(\Omega)$ satisfy $vw = 0$. 
    For every $\theta \in [0,2\pi)$ one has
    \begin{align*}
        0 
        = 
        \mf{a}_{\tilde B}[\modulus{v}, \modulus{w}] 
        \le 
        \Re \big(e^{i\theta} \mf{a}_B [v,w] \big)
        ,
    \end{align*}
    where the equality uses the locality of $\mf{a}_{\tilde B}$ and the inequality follows from the domination assumption between the semigroups and Proposition~\ref{prop:domination-robin}\ref{prop:domination-robin:itm:domination} and~\ref{prop:domination-robin:itm:forms} since $e^{i\theta} v \overline{w} = 0$. 
    Hence, $\mf{a}_B [v,w] = 0$.
\end{proof}

If the operator $L_B$ in Setting~\ref{setting:main} satisfies the estimate~\eqref{eq:neumann-domination}, then Proposition~\ref{prop:local-by-dom}, applied to the operator $\tilde B \coloneqq 0$, shows that $B$ is a multiplication operator, so $L_B$ is subject to local Robin boundary conditions. 
More generally, this argument shows that the problem of Gaussian estimates in the case of non-local Robin boundary conditions cannot be tackled simply using domination by other semigroups associated to local forms. 

Similar results as in Proposition~\ref{prop:local-by-dom} were shown in~\cite[Corollary 4.2]{Akhlil} and \cite[Theorem~3.2]{AroraChillDjida}. 
In those references, the dominated semigroup is assumed to be positive. 
This assumption is not needed in Proposition~\ref{prop:local-by-dom} since the form domains of $\mf{a}_B$ and $\mf{a}_{\tilde B}$ coincide.

\section{Eventual positivity: the case \texorpdfstring{$\rs(-L_B) = 0$}{}}
\label{sec:spectrum-epos}

We turn our attention to the question of eventual positivity of the semigroup $(e^{-tL_B})_{t\ge 0}$. As mentioned in Section~\ref{sec:epos} in the introduction, the key tools to establish this property are ultracontractivity and the existence of positive eigenvector for the spectral bound $\rs(-L_B)$ of the generator $-L_B$.
A more detailed analysis of the spectrum of $-L_B$ is therefore required, and for this purpose we distinguish two cases for the spectral bound of the operator $-L_B$. In this section, we consider $\rs(-L_B)=0$, while Section~\ref{sec:symm} covers the case $\rs(-L_B)>0$. Actually, in both cases we are able to reformulate the spectral bound condition in terms of the spectrum of the boundary operator $B$. This is likely to be more practical, since one expects to have more explicit information about the boundary conditions in concrete examples.

We first prove a result concerning the triviality of the peripheral spectrum of the operator $-L_B$.

\begin{theorem}
	\label{thm:A-spectrum}
	Let the conditions of Setting~\ref{setting:main} be satisfied, and assume that $B+B^*$ is positive semi-definite on $L^2(\partial\Omega)$. The following assertions hold:
    \begin{enumerate}[\upshape (i)]
        \item\label{thm:A-spectrum:it:spec-bound-negative} $\rs(-L_B)\le 0$ and $\sigma(-L_B)\cap i\RR \subseteq \{0\}$.
        
        \item\label{thm:A-spectrum:it:spec-bound-zero} Let $\sigma(-L_B)\cap i\RR\ne\emptyset$. Then $\rs(-L_B) = 0$, $\sigma(-L_B)\cap i\RR = \{0\}$, and every eigenfunction corresponding to the eigenvalue $0$ is constant. Consequently $\dim\ker(-L_B) = 1$.

        \item\label{thm:A-spectrum:it:B} $0 \in \sigma(-L_B)$ if and only if $B\mathbf{1}_{\partial\Omega} = 0$.
    \end{enumerate}
    Moreover, all the above assertions also hold for the adjoint operator $-L_B^*$.
\end{theorem}

\begin{proof}
	\ref{thm:A-spectrum:it:spec-bound-negative}
    The assertion $\rs(-L_B)\le 0$ is a consequence of Proposition~\ref{prop:basic}\ref{prop:basic:itm:contractive}: since $B+B^*$ is positive semi-definite, it follows that the semigroup $(e^{-tL_B})_{t\ge 0}$ is contractive, and then by a standard result in semigroup theory~\cite[Proposition IV.2.2]{EN00} we conclude
	\begin{equation*}
		\rs(-L_B) \le \omega_0(-L_B) \le 0.
	\end{equation*}
	
	\ref{thm:A-spectrum:it:spec-bound-zero}
    We know from Proposition~\ref{prop:basic}\ref{prop:basic:itm:operator}
    that $-L_B$ has compact resolvent, and hence $\sigma(-L_B)$ consists only of eigenvalues. Assume that $-L_B$ has an eigenvector $v\in \dom(-L_B)$ such that $\|v\|_2 =1$ and $-L_B v = i\omega v$ for some $\omega\in\RR$. Then
	\begin{equation*}
		-i\omega = \braket{L_B v, v}_\Omega = \int_\Omega A\nabla v \cdot \overline{\nabla v}\,dx + \braket{B\gamma(v),\gamma(v)}_{\partial\Omega}.
	\end{equation*}
	Therefore, upon taking the real part of the above equation and using the uniform ellipticity~\ref{setting:main:itm:unif-elliptic}, we obtain
	\begin{equation*}
		0\le \alpha \|\nabla v\|^2_2 \le \Re\int_\Omega A\nabla v\cdot\overline{\nabla v}\,dx = -\Re\braket{B\gamma(v),\gamma(v)}_{\partial\Omega} \le 0,
	\end{equation*}
	where the last inequality is due to the positive semi-definiteness of $B + B^*$. Since $\alpha>0$, it follows that $\nabla v = 0$, so $v$ is a non-zero constant function. Consequently
    \[
        0 = \mf{a}_B[v,\varphi] = \braket{L_B v,\varphi}_\Omega
    \]
    for all $\varphi \in C^\infty_c(\Omega) \subseteq H^1(\Omega)$, which shows that $L_B v = 0$ and hence $\omega=0$. Finally, since the semigroup $(e^{-tL_B})_{t\ge 0}$ is analytic (Proposition~\ref{prop:basic}\ref{prop:basic:itm:operator}), it is known that $\rs(-L_B)=\omega_0(-L_B)$, see~\cite[Corollary IV.3.12]{EN00}. However, $0\in\sigma(-L_B)$ also implies $\rs(-L_B)\ge 0$, and hence $\rs(-L_B)=0$.
    
    \ref{thm:A-spectrum:it:B}
    If $0\in\sigma(-L_B)$, then part~\ref{thm:A-spectrum:it:spec-bound-zero} shows that $\rs(-L_B)=0$, and $\mathbf{1}$ spans the associated eigenspace. Then
    \[
        0 = \braket{L_B \mathbf{1},v}_{\Omega} =  \mf{a}_B[\mathbf{1},v] = \braket{B\mathbf{1}_{\partial\Omega}, \gamma(v)}_{\partial\Omega}
    \]
    for all $v\in H^1(\Omega)$. Since the image of $H^1(\Omega)$ under the trace operator $\gamma$ is known to be the fractional Sobolev space $H^{1/2}(\partial\Omega)$ -- see for instance~\cite[pp.\ 98--99 and Theorem 3.37]{McLean} -- and this space is dense in $L^2(\partial\Omega)$, it follows that $B\mathbf{1}_{\partial\Omega} = 0$.
    
    Conversely, if $B\mathbf{1}_{\partial\Omega}=0$, then
	\begin{equation*}
		\mf{a}_B[\mathbf{1},\varphi] = \braket{B\mathbf{1}_{\partial\Omega},\gamma(\varphi)}_{\partial\Omega} = 0 = \braket{0,\varphi}_{\Omega}
	\end{equation*}
	for all $\varphi\in H^1(\Omega)$. This shows that $\mathbf{1}\in\dom(L_B)$ with $L_B \mathbf{1} = 0$, and hence $0\in\sigma(-L_B)$.

    Finally, assertions~\ref{thm:A-spectrum:it:spec-bound-negative}--\ref{thm:A-spectrum:it:B} also hold for $-L_B^*$ (recall Proposition~\ref{prop:basic}\ref{prop:basic:itm:adjoint} for the precise description of the adjoint generator), because clearly the positive semi-definiteness of $B+B^*$ is unchanged if $B$ is replaced by $B^*$.
\end{proof}

In the one-dimensional case, where $\Omega=(a,b)$ is a bounded open interval, $L^2(\partial\Omega)$ can be identified with $\CC^2$. We then obtain a perhaps surprising corollary for positivity of the semigroup.

\begin{corollary}
	Assume that all the conditions of Setting~\ref{setting:main} are satisfied and let the dimension be $d=1$. 
    If $B+B^*$ is positive semi-definite (on $\CC^2$) and $(1\quad 1)^\top\in\ker B$, then $(e^{-tL_B})_{t\ge 0}$ is positive.
\end{corollary}

\begin{proof}
	By Theorem~\ref{thm:A-spectrum}\ref{thm:A-spectrum:it:B}, the condition $(1\quad 1)^\top \in\ker B$ is equivalent to $0\in\sigma(-L_B)$. 
    Hence part~\ref{thm:A-spectrum:it:spec-bound-zero} of that theorem shows that $\sigma(-L_B)\cap i\RR = \{0\}$, $\rs(-L_B)=0$, and any eigenvector $v$ corresponding to $0$ is a non-zero constant function on $\Omega$. 
    The positive semi-definiteness of $B+B^*$ implies that the diagonal entries of $B$ are real and non-negative, 
    and the condition $(1\quad 1)^\top \in\ker B$ then implies that the off-diagonal entries of $B$ are real and non-positive. 
    In other words, the matrix $-B$ is real and has non-negative off-diagonal entries. 
    Therefore we can choose $c\in\RR$ such that $-B + cI \ge 0$, and thus 
    $
        e^{-tB} 
        = 
        e^{-tc}e^{t(-B+cI)} 
        \ge 
        0
    $
    for every $t \ge 0$, where the inequality follows from the series representation of the matrix exponential function. 
    Hence $(e^{-tL_B})_{t\ge 0}$ is a positive semigroup by Proposition~\ref{prop:positivity}.
\end{proof}

In higher dimensions we cannot expect positivity, but only uniform eventual positivity of the semigroup $(e^{-tL_B})_{t\ge 0}$. 
Proving this property is possible due our earlier result on ultracontractivity (Theorem~\ref{thm:ultra}), the spectral analysis in Theorem~\ref{thm:A-spectrum}, and general criteria for eventual positivity of semigroups from the paper~\cite{DG18}.
To keep the presentation self-contained, we explicitly state the following Proposition~\ref{prop:unif-ev-pos}, which is a special case of~\cite[Theorem~3.1]{DG18} suitable for the current purposes. 
Note that in the assumptions of \cite[Theorem~3.1]{DG18} dual operators are understood in the sense of Banach space duals, while the assumptions of Proposition~\ref{prop:unif-ev-pos} below are stated in terms of Hilbert space adjoints of operators. 
In the setting of Proposition~\ref{prop:unif-ev-pos} those concepts coincide since all involved operators are real. 
We refer to the paragraph before \cite[Corollary~3.5]{DG18} for a detailed discussion of this.

In condition~\ref{prop:unif-ev-pos:itm:spectral} in the following proposition, the spectral bound of the operator $C$ is said to be a \emph{dominant spectral value} if $\rs(-C)\in\sigma(-C)$ and $\sigma(-C) \cap (\rs(-C) + i\RR) = \{ \rs(-C) \}$.
We note that condition~\ref{prop:unif-ev-pos:itm:smoothing} in the proposition implies that the operators $e^{-tC}$ on $L^2(M,\mu)$ are compact for all sufficiently large $t$ (see e.g.\ \cite[Corollary~2.4]{DG18}). 
Hence, all spectral values of $C$ are eigenvalues (see e.g.\ \cite[Corollary~V.3.2 and assertion~(ii) before Proposition~IV.1.18]{EN00}).

\begin{proposition}
    \label{prop:unif-ev-pos}
    Let $(M,\mu)$ be a finite measure space, and let $(e^{-tC})_{t\ge 0}$ be a real $C_0$-semigroup on $L^2(M,\mu)$ such that the generator $-C$ has non-empty spectrum. Assume that the following conditions hold:
    \begin{enumerate}[\upshape (i)]
        \item\label{prop:unif-ev-pos:itm:smoothing} 
        There exist times $t_1, t_2 \ge 0$ such that $e^{-t_1 C}(L^2(M,\mu)) \subseteq L^\infty(M,\mu)$ and $e^{-t_2 C^*}(L^2(M,\mu)) \subseteq L^\infty(M,\mu)$.
        
        \item\label{prop:unif-ev-pos:itm:spectral} 
        The spectral bound $\rs(-C)$ is a dominant spectral value, the corresponding eigenspace $\ker(\rs(-C) + C)$ is one-dimensional and spanned by an eigenfunction $u$ such that $u \ge \mathbf{1}$, and the dual eigenspace $\ker(\rs(-C) + C^*)$ contains an eigenfunction $v$ such that $v \ge \mathbf{1}$.
    \end{enumerate}
    Then there exist $t_0 \ge 0$ and $\delta > 0$ such that
    \[
        e^{-t\rs(-C)} e^{-tC}f \ge \delta \left(\int_M f\,d\mu \right) \mathbf{1}
    \]
    for all $t\ge t_0$ and for all $0\le f\in L^2(M,\mu)$. 
    In particular, the semigroup $(e^{-tC})_{t\ge 0}$ is uniformly eventually positive.
\end{proposition}

Note that the factor $e^{-t\rs(-C)}$ on the left hand side of the conclusion of Proposition~\ref{prop:unif-ev-pos} does not occur in \cite[Theorem~3.1]{DG18}. 
This is a typo in this reference. 
It is easy to see that, without this factor, the estimate does not scale correctly if one adds a multiple of the identity to $C$.

We can now present a sufficient criterion for uniform eventual positivity of the semigroup $(e^{-tL_B})_{t\ge 0}$.

\begin{theorem}
	\label{thm:eventual-pos}
	In Setting~\ref{setting:main}, assume that the dimension is $d\ge 2$.
    Assume that $B$ is real and acts boundedly on $L^1(\partial\Omega)$ and $L^\infty(\partial\Omega)$. 
    Assume moreover that $B+B^*$ is positive semi-definite and $B\mathbf{1}_{\partial\Omega}=0$. Then there exist $t_0\ge 0$ and a constant $\delta>0$ such that
	\begin{equation*}
		e^{-tL_B}f \ge \delta \left(\int_\Omega f\,dx\right)\mathbf{1} \qquad\forall\, t\ge t_0
	\end{equation*}
	for all $0\le f\in L^2(\Omega)$. 
    In particular, $(e^{-tL_B})_{t \ge 0}$ is uniformly eventually positive.
\end{theorem}

\begin{proof}
    We check the assumptions of Proposition~\ref{prop:unif-ev-pos}. 
    Due to assumption~\ref{setting:main:itm:A-real} in Setting~\ref{setting:main} and since $B$ is real, it follows from Proposition~\ref{prop:basic}\ref{prop:basic:itm:real} that the semigroup $(e^{-tL_B})_{t\ge 0}$ is real.
    Since $B$ acts boundedly on $L^1$ and $L^\infty$, it follows from Theorem~\ref{thm:ultra} that the ultracontracitvity condition~\ref{prop:unif-ev-pos:itm:smoothing} in Proposition~\ref{prop:unif-ev-pos} holds. 
    As $B\mathbf{1}_{\partial\Omega} = 0$, we conclude from Theorem~\ref{thm:A-spectrum} that $\rs(-C) = 0$ and that the spectral condition~\ref{prop:unif-ev-pos:itm:spectral} in Proposition~\ref{prop:unif-ev-pos} also holds.  
    So Proposition~\ref{prop:unif-ev-pos} can be applied and gives the claimed estimate.
\end{proof}

We now give some explicit examples of operators $B \in \mathcal{L}(L^2(\Omega))$ that satisfy the assumptions of Theorem~\ref{thm:eventual-pos} but for which the semigroup $(e^{-tB})_{t\ge 0}$ is not positive.

\begin{examples}
    \label{exa:operators-B-ev-pos}
    \begin{enumerate}[(a)]
    	\item 
        Let $0\ne v\in L^\infty(\partial\Omega)$ be a real-valued function satisfying $\int_{\partial\Omega} v\,d\sigma = 0$, and consider the rank-$1$ operator $B \in \mathcal{L}(L^2(\partial \Omega))$ defined by
    	\begin{equation*}
    		Bf \coloneqq (v\otimes v)f = \left(\int_{\partial\Omega} fv \,d\sigma\right) v \qquad \forall\,f\in L^2(\partial\Omega).
    	\end{equation*}
    	This is clearly a bounded operator on $L^2(\partial\Omega)$, and acts boundedly on $L^1$ and $L^\infty$. 
        By construction it holds that $B\mathbf{1}_{\partial\Omega} = 0$. 
        As $v$ is real, $B$ is self-adjoint and real, and moreover
    	\begin{equation*}
    		\braket{(v\otimes v)f,f}_{\partial\Omega} = \left(\int_{\partial\Omega}vf\,d\sigma\right) \left(\int_{\partial\Omega} v\overline{f}\,d\sigma\right) = \left\vert \int_{\partial\Omega} vf\,d\sigma\right\vert^2 \ge 0
    	\end{equation*}
    	for all $f\in L^2(\partial\Omega)$. 
        Thus $B+B^*=2B$ is positive semi-definite. 
        To show that $-B$ does not generate a positive semigroup, we test the form $\mf{b}$ from Setting~\ref{setting:main} with $v$ itself and obtain
    	\begin{equation*}
    		\mf{b}[v^+, v^-] = \braket{Bv^+, v^-}_{\partial\Omega} = \left(\int_{\partial\Omega} \modulus{v^+}^2 \,d\sigma\right) \left(\int_{\partial\Omega} \modulus{v^-}^2 \,d\sigma\right) > 0,
    	\end{equation*}
    	which shows that the form $\mf{b}$ violates the Beurling--Deny criterion (see \cite[Theorem 2.6]{Ouh}). 
        Hence $(e^{-tB})_{t\ge 0}$ is not positive.
    	
    	\item 
        Let us also give an example where $B$ is not a kernel operator. Consider the unit disk $\Omega = \{x\in\RR^2 : \modulus{x}<1\}$ with boundary $\partial\Omega = \{x\in\RR^2 : \modulus{x}=1\}$, and let $R \in \RR^{2 \times 2}$ be the matix that rotates all elements of $\RR^2$ anticlockwise by the angle $\frac{\pi}{2}$.
        
        Define the operator $B\in\mathcal{L}(L^2(\partial\Omega))$ by $Bf \coloneqq f\circ R - f\circ R^*$.
        Then $B+B^* = 0$ and hence $B$ is positive semi-definite, albeit in a trivial way. In addition, it is clear that $B$ is real, acts boundedly on $L^1(\partial\Omega)$ and $L^\infty(\partial\Omega)$, and $B\mathbf{1}_{\partial\Omega} = 0$. 
        
        Again, let us show that $(e^{-tB})_{t\ge 0}$ is not positive using the Beurling--Deny criterion. 
        Let $f\in L^2(\partial\Omega)$ be the function that is $1$ in the first quadrant, $-1$ in the fourth quadrant, and $0$ everywhere else.
        A simple calculation yields $\mf{b}[f^+, f^-] = \braket{Bf^+, f^-} = \frac{\pi}{2} > 0$, hence the form $\mf{b}$ does not satisfy the Beurling--Deny criterion.
    \end{enumerate}
\end{examples}

Note that for the operators $B$ in Example~\ref{exa:operators-B-ev-pos}, ultracontractivity of $(e^{-tL_B})_{t \ge 0}$ --- which is essential to get eventual positivity --- cannot be shown by means of the domination result in \cite[Theorem~3.3]{GeMiNiOu}. 
Indeed, for this one would need a negative operator on $L^2(\partial\Omega)$, called $\Theta_1$ in \cite{GeMiNiOu}, such that
\begin{align*}
    \modulus{e^{-tL_B} f} \le e^{-tL_{\Theta_1}}\modulus{f} \le e^{t\Delta_N} \modulus{f}
\end{align*}
for all $f \in L^2(\Omega)$ and $t \ge 0$, where $\Delta_N$ denotes the Neumann Laplacian. 
As explained after Proposition~\ref{prop:local-by-dom}, this implies that $B$ is a multiplication operator.
However, this is clearly not the case for the operators $B$ in Example~\ref{exa:operators-B-ev-pos}.

\section{Eventual positivity: the case \texorpdfstring{$\rs(-L_B)>0$}{}}
\label{sec:symm}

In Theorem~\ref{thm:A-spectrum}, the spectral condition $\rs(-L_B)=0$ conveniently led to a positive constant eigenfunction, which in turn yielded the uniform eventual positivity in Theorem~\ref{thm:eventual-pos}, but the arguments cannot be adapted to the case $\rs(-L_B)>0$. Thus, instead of developing a general theory, we change the perspective of our analysis in this section: we show how symmetry conditions on the domain $\Omega$ and the coefficients of the differential operator can yield a positive leading eigenfunction of the operator $-L_B$. The interaction between symmetry and spectral theory is a classical area of study. Indeed, quite general properties of eigenfunctions of the Dirichlet Laplacian on symmetric domains were extensively investigated by Pereira~\cite{Per}.

\subsection{Invariance and equivariance under group actions}
\label{sec:invariance-group-actions}

Let $O(d) \subseteq \RR^{d \times d}$ denote the group of orthogonal $d \times d$ matrices and let $G \subseteq O(d)$ be a closed subgroup. 
We say that a subset $S \subseteq \RR^d$ is \emph{$G$-invariant} if $g(S) \subseteq S$ for all $g \in G$. 
As $G$ is a group, this is equivalent to the formally stronger condition $g(S) = S$ for all $g \in G$.

Now let $\emptyset \not= \Omega \subseteq \RR^d$ be open and $G$-invariant. 
Then the boundary $\partial \Omega$ is $G$-invariant as well, since every $g \in G$ is a homeomorphism on $\RR^d$. 
As every $g$ preserves the Lebesgue measure, the linear operator $T_g \colon L^2(\Omega) \to L^2(\Omega)$ given by 
\begin{align*}
    (T_g u)(x) \coloneqq u(g^{-1} x)
\end{align*}
for all $u \in L^2(\Omega)$ and for $x \in \Omega$ is well-defined and unitary. 
In fact, it is easy to verify that the map $G\ni g\mapsto T_g \in\mathcal{L}(L^2(\Omega))$ is a representation of $G$, and this is the reason for using $g^{-1}$ instead of $g$ in the definition. 
In Pereira's terminology~\cite{Per}, this is called the \emph{quasi-regular representation} of $G$. 
 
The chain rule for Sobolev functions (see e.g.~\cite[Proposition 9.5]{Bre}) yields $\nabla (u\circ g^{-1}) = g (\nabla u \circ g^{-1})$ for all $g\in G$ and $u\in H^1(\Omega)$. Using this, it easy to check that each operator $T_g$ leaves $H^1(\Omega)$ invariant and is isometric with respect to the $H^1(\Omega)$ norm.

The fact that each $g \in G$ acts isometrically on $\RR^d$ also implies that it preserves the $(d-1)$-dimensional Hausdorff measure on $\partial \Omega$. 
Hence, for each $g\in G$, the linear operator $S_g \colon L^2(\partial \Omega) \to L^2(\partial \Omega)$, given by
\begin{align*}
    (S_g f)(x) \coloneqq f(g^{-1} x)
\end{align*}
for all $f \in L^2(\partial\Omega)$ and for $x \in \partial \Omega$, is well-defined and unitary. 
If $\Omega$ is a bounded Lipschitz domain, the relation between $T_g$ and $S_g$ in terms of the trace operator is described in assertions~\ref{prop:symmetry:itm:bdry-intertwine}--\ref{prop:symmetry:itm:trace} of Proposition~\ref{prop:symmetry} in the appendix.

Symmetry properties of functions and operators with respect to the action of the group $G$ can be described by the notions of \emph{$G$-invariance} and \emph{$G$-equivariance}, which are defined as follows. 
We say that a function $u \in L^2(\Omega)$ is \emph{$G$-invariant} if $T_g u = u$ for all $g \in G$. 
We use the notation 
\begin{equation}
    \label{eq:symm-subspace}
	F_G 
    \coloneqq 
    \{ u\in L^2(\Omega) : T_g u = u \quad\forall\,g\in G\}
\end{equation}
for the fixed space of the operator family $(T_g)_{g \in G}$ on $L^2(\Omega)$. 
In other words, $F_G$ consists of all $G$-invariant functions in $L^2(\Omega)$. 
Clearly, $F_G$ is a closed vector subspace of $L^2(\Omega)$. 
A number of useful properties of $F_G$, its orthogonal complement $F_G^\perp$, and their relation to the Sobolev space $H^1(\Omega)$ are proved in Proposition~\ref{prop:symmetry} in the appendix.

A bounded linear operator $B \in \mathcal{L}(L^2(\partial\Omega))$ is called \emph{$G$-equivariant} if 
\begin{align*}
    B S_g f = S_g Bf \qquad \, \forall f \in L^2(\partial \Omega).
\end{align*}
A linear operator $C \colon \dom(C)\subseteq L^2(\Omega) \to L^2(\Omega)$ is called \emph{$G$-equivariant} if for all $u\in \dom(C)$ and $g\in G$ it holds that
\begin{equation}
	\label{eq:G-equi}
	T_g u\in \dom(C) \quad\text{and}\quad CT_g u = T_g Cu.
\end{equation}
In regards to spectral theory, the following observation is crucial. If $C\colon \dom(C)\subseteq L^2(\Omega) \to L^2(\Omega)$ is a $G$-equivariant operator and $u$ is an eigenfunction corresponding to an eigenvalue $\lambda$, then we observe
\begin{equation*}
	C T_g u = T_g Cu = \lambda T_g u \qquad\forall\, g\in G,
\end{equation*}
which shows that $T_g u$ is also an eigenfunction corresponding to $\lambda$. 
Hence the eigenspace $\ker(\lambda - C)$ is invariant under all the operators $T_g$. 

A sesquilinear form $\mf{c} \colon \dom(\mf{c})\times \dom(\mf{c})\subseteq L^2(\Omega) \times L^2(\Omega)\to \CC$ is called \emph{$G$-invariant} if for all $u, v \in \dom(\mf{c})$ and $g\in G$, it holds that
\begin{equation}
    \label{def:G-equi-form:eq}
    T_g u\in \dom(\mf{c}) \quad\text{and}\quad \mf{c}[T_g u, T_g v] = \mf{c}[u, v].
\end{equation}

\begin{proposition}
    \label{prop:G-equi-op-form}
    Let $G$ be a closed subgroup of the orthogonal group $O(d)$, and suppose $\mf{c} \colon \dom(\mf{c})\times \dom(\mf{c})\subseteq L^2(\Omega) \times L^2(\Omega) \to \CC$ is a densely-defined sesquilinear form. 
    If $\mf{c}$ is $G$-invariant, then its associated operator $C \colon \dom(C)\subseteq L^2(\Omega) \to L^2(\Omega)$ is $G$-equivariant in the sense of~\eqref{eq:G-equi}.
\end{proposition}

\begin{proof}
    Suppose $u\in\dom(C)$, and let $f \coloneqq Cu$. 
    Since $\mf{c}$ is $G$-invariant, we obtain
    \[
        \mf{c}[T_g u, T_g v] = \mf{c}[u,v] = \braket{f,v} = \braket{T_g f, T_g v}
    \]
    for all $v\in\dom(\mf{c})$, where we have used that $g \in O(d)$ in the last equality. 
    Since we can write $w = T_g (T_{g^{-1}}w)$ for every $w\in \dom(\mf{c})$, it follows that $\mf{c}[T_g u, w] = \braket{T_g f, w}$ for all $w\in \dom(\mf{c})$ and $g\in G$. Hence $T_g u \in \dom(C)$ with $CT_g u = T_g f = T_g Cu$, which shows that $C$ is $G$-equivariant.
\end{proof}

\subsection{Equivariance of the Robin Laplacian under group actions}

We turn our attention to the analysis of the operator $L_B$ from Setting~\ref{setting:main} when the domain $\Omega \subseteq \RR^d$ is invariant under a closed subgroup $G \subseteq O(d)$. 
Our assumptions here are more specific than in the previous sections, so we highlight them separately.

\begin{setting}
    \label{setting:symmetry}
    Let all assumptions of Setting~\ref{setting:main} be satisfied and let $G$ be a closed subgroup of the orthogonal group $O(d)$. 
    Assume that the domain $\Omega$ is $G$-invariant and that the following properties hold.
    \begin{enumerate}[label=(B\arabic*)]
        \item\label{setting:symmetry:itm:coefficients} 
        For each $x \in \Omega$ the matrix $A(x) \in \RR^{d \times d}$ is symmetric (i.e.\ $A(x)^\top = A(x)$) and commutes with each $g\in G$ (i.e.\ $gA(x) = A(x)g$). 
        Moreover, $A \colon \Omega \to \RR^{d \times d}$ is $G$-invariant, i.e.\ $A(gx) = A(x)$ for all $x \in \Omega$ and $g \in G$.
        
        \item\label{setting:symmetry:itm:B2} 
        $B$ is a self-adjoint and $G$-equivariant operator.
    \end{enumerate}
\end{setting}

Assumption~\ref{setting:symmetry:itm:coefficients} in Setting~\ref{setting:symmetry} is quite strong if the group $G$ is sufficiently large. 
To explain this, we mention some observations from representation theory in the following proposition; for the convenience of the reader, we include the simple proofs. Recall that an action of a group $G$ on a set $X$ is called \emph{transitive} if for every pair of elements $x,y\in X$, there exists $g\in G$ such that $g\cdot x=y$.

\begin{proposition}
    \label{prop:invariant-coeff-trivial}
    Let $G \subseteq O(d)$ be a group and assume that $A \in \RR^{d \times d}$ commutes with all $g \in G$.
    \begin{enumerate}[\upshape(i)]
        \item\label{prop:invariant-coeff-trivial:itm:full} 
        If $G = O(d)$, then $A$ is a scalar multiple of the identity matrix $\id$.

        \item\label{prop:invariant-coeff-trivial:itm:transitive} 
        If $G$ acts transitively on the Euclidean unit sphere of $\RR^d$ and $A$ is symmetric, then $A$ is a scalar multiple of $\id$.
    \end{enumerate}
\end{proposition}

\begin{proof}
    \ref{prop:invariant-coeff-trivial:itm:full} 
    Let $\mathcal{C}$ denote the sub-algebra of $\RR^{d \times d}$ generated by $O(d)$. 
    As $A$ commutes with all elements of $\mathcal{C}$, it suffices to show that $\mathcal{C} = \RR^{d \times d}$. 
    If $C \in \RR^{d \times d}$ is skew-symmetric, one has $\exp(tC) \in O(d)$ for all $t \in \RR$ thus, $C = \frac{d \exp(tk)}{dt}\big|_{t = 0} \in \mathcal{C}$. 
    Now suppose that $C \in \RR^{d \times d}$ is symmetric. 
    After rescaling, we may assume that the eigenvalues of $C$ (which are all real) lie in the interval $[-1,1]$. By the spectral theorem for symmetric matrices, there exist $g\in O(d)$ and a diagonal matrix $D \in \RR^{d \times d}$ with entries in $[-1,1]$ such that $C=g D g^\top$. 
    The matrix $D$ can then be expressed as a convex combination of diagonal matrices with entries in $\{-1,1\}$, and thus $C\in\mathcal{C}$.
    Since every $C \in \RR^{d \times d}$ can be written as the sum of a symmetric and a skew-symmetric matrix, it follows that $\mathcal{C} = \RR^{d \times d}$, as claimed.
    
    \ref{prop:invariant-coeff-trivial:itm:transitive} 
    Since $A \in \RR^{d \times d}$ is symmetric, it is diagonalisable and has only real eigenvalues. 
    It suffices to prove that all eigenvalues of $A$ are equal. 
    So let $\lambda,\mu \in \RR$ be eigenvalues of $A$ 
    with normalised eigenvectors $x,y \in \RR^d$ for $\lambda$ and $\mu$, respectively. 
    Since $G$ acts transitively on the unit sphere, there exists $g \in G$ such that $gx = y$. 
    Hence 
    \begin{align*}
        \mu y 
        = 
        Ay 
        = 
        A gx 
        = 
        g Ax 
        = 
        g \lambda x 
        = 
        \lambda y
        ,
    \end{align*}
    which shows that $\lambda = \mu$, as claimed.
\end{proof}

Under the conditions of Setting~\ref{setting:symmetry}, it is straightforward to verify that the form $\mf{a}_B$ is symmetric and $G$-invariant. 
Hence $L_B$ is self-adjoint, and by Proposition~\ref{prop:G-equi-op-form}, also $G$-equivariant. 
In contrast to the results in Section~\ref{sec:spectrum-epos}, we no longer require the assumption that $B+B^*$ is positive semi-definite, but we require $L_B$ to be self-adjoint in order to employ a variational principle in the following result.

\begin{theorem}
	\label{thm:symm-eig}
	Let the conditions of Setting~\ref{setting:symmetry} be satisfied, and write $L^2_0(\partial\Omega)$ for the subspace of $L^2(\partial\Omega)$ consisting of mean-zero functions, i.e.\ $f\in L^2(\partial\Omega)$ such that $\int_{\partial\Omega}f\,d\sigma=0$. 
    Assume that
	\begin{equation}
		\label{eq:sB-positive} 
        \braket{B\mathbf{1}_{\partial\Omega},\mathbf{1}_{\partial\Omega}}_{L^2(\partial\Omega)} < 0.
	\end{equation}
	Then $\rs(-L_B)>0$ and there exists a constant $\beta>0$ (depending only on the domain $\Omega$ and the uniform ellipticity constant $\alpha$ in~\ref{setting:main:itm:unif-elliptic}) with the following property: 
    if
    \begin{equation*}
        \norm{B|_{L^2_0(\partial\Omega)} }_{L^2_0(\partial\Omega) \to L^2(\partial\Omega)}\le \beta,
    \end{equation*}
    then every eigenfunction of $-L_B$ associated to $\rs(-L_B)$ is $G$-invariant.
\end{theorem}

\begin{proof}
	As in~\eqref{eq:symm-subspace}, we write $F_G \subseteq L^2(\Omega)$ for the space of all $G$-invariant functions in $L^2(\Omega)$. 
    Clearly $\mathbf{1}\in F_G\cap H^1(\Omega)$, so if $u \in F_G^\perp \cap H^1(\Omega)$, then $\int_\Omega u \, dx = 0$ and hence the Poincar\'{e} inequality~\cite[Theorem 13.27]{Leoni} implies
	\begin{equation*}
		\|u\|^2_2 \le c_0 \|\nabla u\|^2_2 \qquad\forall\, u\in F_G^\perp\cap H^1(\Omega)
	\end{equation*}
	with a constant $c_0>0$ depending only on $\Omega$. We combine this result with the standard trace inequality $\|\gamma(u)\|^2_{2,\partial\Omega} \le c_1 \|u\|^2_{H^1(\Omega)}$ to obtain
	\begin{equation*}
		\|\gamma(u)\|^2_{2,\partial\Omega} \le c_1\|u\|^2_{H^1(\Omega)} \le (c_0+1)c_1\|\nabla u\|^2_2
	\end{equation*}
	for all $u\in F_G^\perp\cap H^1(\Omega)$. 
    Now we set $C\coloneqq (c_0+1)c_1$ and $\beta \coloneqq \alpha C^{-1}$, where $\alpha > 0$ is the ellipticity constant from~\ref{setting:main:itm:unif-elliptic}. 
    Then 
	\begin{align}
        \label{thm:symm-eig:eq:aB-lower-bound}
		\mf{a}_B[u,u] &\ge \alpha\|\nabla u\|^2_2 + \braket{B\gamma(u),\gamma(u)}_{\partial\Omega} \notag \\
        &\ge \left(\alpha C^{-1}\|\gamma(u)\|_{2,\partial\Omega} - \|B\gamma(u)\|_{2,\partial\Omega} \right) \|\gamma(u)\|_{2,\partial\Omega}
	\end{align}
	for all $u\in F_G^\perp\cap H^1(\Omega)$. One can show that $\gamma(u)\in L^2_0(\partial\Omega)$ -- this follows from Proposition~\ref{prop:symmetry}\ref{prop:symmetry:itm:trace} in the appendix. Using this property together with the assumption $\norm{B|_{L^2_0(\partial\Omega)} }_{L^2_0(\partial\Omega) \to L^2(\partial\Omega)}\le \beta =\alpha C^{-1}$ and inequality~\eqref{thm:symm-eig:eq:aB-lower-bound}, we find that
	\begin{equation}
		\label{eq:anti-sym-a}
		\mf{a}_B[u,u]\ge 0 \qquad\forall\, u\in F_G^\perp\cap H^1(\Omega).
	\end{equation}
	
	Since $B$ is self-adjoint, the form $\mf{a}_B$ is symmetric. Thus, by a standard variational principle, the smallest eigenvalue of $L_B$ is given by the Rayleigh quotient
	\begin{equation*}
		-\rs(-L_B) = \min_{0\ne u\in H^1(\Omega)} \frac{\mf{a}_B[u,u]}{\|u\|^2_2}.
	\end{equation*}
	Condition~\eqref{eq:sB-positive} implies $\rs(-L_B)>0$. Indeed, if we use the constant function $\mathbf{1}$ in the form $\mf{a}_B$, we obtain
	\begin{equation*}
		-\rs(-L_B) \le \frac{\mf{a}_B[\mathbf{1},\mathbf{1}]}{\|\mathbf{1}\|^2_2} = \frac{1}{\modulus{\Omega}} \braket{B\gamma(\mathbf{1})),\gamma(\mathbf{1})}_{\partial\Omega} < 0.
	\end{equation*}
    We complete the proof by considering a normalised eigenfunction $u$ of $-L_B$ associated to $\rs(-L_B)$ and showing that $u \in F_G$. 
    Note that $u\in H^1(\Omega)$ is a minimiser of the Rayleigh quotient. 
    Then we can write $u = u_0 + u_1$ where $u_0 \in F_G \cap H^1(\Omega)$ and $u_1 \in F_G^\perp \cap H^1(\Omega)$, see Proposition~\ref{prop:symmetry}\ref{prop:symmetry:itm:H1-directsum} in the appendix. 
    Thus
    \[
        0 
        > 
        -\rs(-L_B) 
        = 
        \mf{a}_B[u,u] 
        = 
        \mf{a}_B[u_0,u_0] + \mf{a}_B[u_1,u_1] 
        \ge 
        \mf{a}_B[u_0,u_0]
        ,
    \]
    where the second equality is due to Proposition~\ref{prop:symmetry}\ref{prop:symmetry:itm:form} and the inequality at the end follows from~\eqref{eq:anti-sym-a}. 
    If $u_1 \ne 0$, then $1 = \norm{u} > \norm{u_0}$ and hence
    \[
        -\rs(-L_B) 
        \ge 
        \mf{a}_B[u_0,u_0] 
        > 
        \frac{\mf{a}_B[u_0,u_0]}{\|u_0\|^2_2}
        .
    \]
    This contradicts that $u$ is a minimiser, so we deduce that $u_1 = 0$, and therefore $u = u_0$ is an element of $F_G$.
\end{proof}

Observe that the last part of the proof of Theorem~\ref{thm:symm-eig} does not really use that we work with the eigenvalue $\rs(-L_B)$. 
The same argument show that

If $\Omega$ is a ball and the group $G$ is sufficiently large, we can say much more about the eigenfunctions of $L_B$ for $\rs(-L_B)$. 
We also obtain sufficient conditions for uniform eventual positivity of the semigroup $(e^{-tL_B})_{t\ge 0}$ in a situation different from Section~\ref{sec:spectrum-epos}.

\begin{theorem}
    \label{thm:ball}
    Let the conditions of Setting~\ref{setting:symmetry} be satisfied, and again write $L^2_0(\partial\Omega)$ for the subspace of $L^2(\partial\Omega)$ consisting of mean-zero functions. 
    Assume in addition that $\Omega$ is the open unit ball in $\RR^d$, that $G$ acts transitively on the unit sphere $\partial\Omega = \mathbb{S}^{d-1}$, 
    and that $B$ satisfies $B(L^\infty(\partial\Omega))\subseteq L^\infty(\partial\Omega)$ and $\braket{B\mathbf{1}_{\partial\Omega},\mathbf{1}_{\partial\Omega}}_{L^2(\partial\Omega)} < 0$. Then there exists a constant $\beta > 0$ (depending only on the domain $\Omega$ and the uniform ellipticity constant $\alpha$ in~\ref{setting:main:itm:unif-elliptic}) such that if
	\begin{equation}
        \label{eq:small-B}
        \norm{B|_{L^2_0(\partial\Omega)} }_{L^2_0(\partial\Omega) \to L^2(\partial\Omega)}\le \beta,
	\end{equation}
	then the following assertions hold:
    \begin{enumerate}[\upshape (i)]
        \item\label{thm:ball:itm:spectrum}
        The eigenvalue $\rs(-L_B)$ of $-L_B$ admits a $G$-invariant eigenfunction $\varphi\in C(\overline{\Omega})$ such that $\varphi(x)\ge c$ for some $c > 0$ and all $x\in\overline{\Omega}$. 
        Moreover, $\rs(-L_B)$ is a simple eigenvalue of $-L_B$ and the only eigenvalue of $-L_B$ in $(0,\infty)$.

        \item\label{thm:ball:itm:ev-pos}
        Assume now that $B\in\mathcal{L}(L^2(\partial\Omega))$ is real. Then there exist $t_0\ge 0$ and a constant $\delta>0$ such that
        \begin{equation*}
    	       e^{-\rs(-L_B)t} e^{-tL_B}f \ge \delta \left(\int_\Omega f \,dx \right) \mathbf{1} \qquad\forall\, t\ge t_0
        \end{equation*}
        for all $0\le f\in L^2(\Omega)$. 
        In particular, the semigroup $(e^{-tL_B})_{t \ge 0}$ is uniformly eventually positive.
    \end{enumerate}
\end{theorem}

Before giving the proof, we briefly comment on the strength of the assumptions of the above theorem.
Since $G$ acts transitively on the unit sphere, it follows from condition~\ref{setting:symmetry:itm:coefficients} and Proposition~\ref{prop:invariant-coeff-trivial}\ref{prop:invariant-coeff-trivial:itm:transitive} that the coefficient matrices satisfy $A(x) = \lambda(x) \id$ for numbers $\lambda(x) \in \RR$ and all $x \in \Omega$. 
Moreover, the invariance condition in~\ref{setting:symmetry:itm:coefficients} implies that $\lambda$ is radially symmetric. 

\begin{proof}[Proof of Theorem~\ref{thm:ball}]
    Choose $\beta > 0$ as in Theorem~\ref{thm:symm-eig}. 
    
    \ref{thm:ball:itm:spectrum}
    Since $B$ is self-adjoint, it follows from the assumption $B(L^\infty(\partial\Omega))\subseteq L^\infty(\partial\Omega)$ and Lemma~\ref{lem:op-L_1-L_infty}\ref{lem:op-L_1-L_infty:itm:adjoint-2} that $B$ acts boundedly on $L^1$ and $L^\infty$.
    Set $\lambda_0 \coloneqq \rs(-L_B)$, and let $\varphi$ be an eigenfunction corresponding to $\lambda_0$. Note that $\lambda_0 > 0$ by Theorem~\ref{thm:symm-eig}.
    We first show that $\varphi\in C(\overline{\Omega})$. Indeed, since $\varphi = e^{-2\lambda_0}e^{-2L_B}\varphi$, we have $\varphi\in L^\infty(\Omega)$ by the ultracontractivity from Theorem~\ref{thm:ultra}.
    Now set $\psi \coloneqq L_B \varphi = -\lambda_0\varphi \in L^\infty(\Omega)$. 
    Moreover, $\varphi$ satisfies
    \[
        \int_\Omega A\nabla \varphi \cdot \overline{\nabla v} \,dx 
        = 
        \mf{a}_B[\varphi,v] 
        - 
        \mf{b}[\gamma(\varphi), \gamma(v)]
        =
        \int_\Omega \psi\overline{v}\,dx - \int_{\partial\Omega} B\gamma(\varphi)\overline{\gamma(v)}\,d\sigma
    \]
    for all $v\in H^1(\Omega)$. The result $\varphi\in C(\overline{\Omega})$ now follows from~\cite[Proposition 3.6]{N11}.
    
    Next, we show that $\varphi$ can be chosen so that $\varphi(x) > 0$ for all $x\in\Omega$.
    According to Theorem~\ref{thm:symm-eig}, the function $\varphi$ is $G$-invariant, 
    so the continuity of $\varphi$ implies that $\varphi(gx) = \varphi(x)$ for all $g \in G$ and all $x \in \overline{\Omega}$. 
    As $G$ acts transitively on the boundary $\mathbb{S}^{d-1}$ of $\Omega$, it follows that $\varphi|_{\partial\Omega}$ is constant. 
    If this constant is $0$, then the properties $\varphi\in C(\overline{\Omega}) \cap H^1(\Omega)$ and $\varphi|_{\partial\Omega} = 0$ imply $\varphi\in H^1_0(\Omega)$.
    Since the equation $(-L_B-\lambda_0)\varphi = 0$ holds with $-\lambda_0 < 0$, it is a consequence of the maximum principle for divergence-form elliptic operators~\cite[Corollary 8.2]{GT} that $\varphi \equiv 0$, which is impossible. Hence, after multiplying by a suitable scalar, we may assume $\varphi|_{\partial\Omega}=1$. Now the maximum principle~\cite[Theorem 8.1]{GT} implies that $\varphi\ge 0$ in $\Omega$.
	Furthermore, the argument in~\cite[Theorem 3.1]{Dan09}, which employs the weak Harnack inequality, shows that $\varphi(x)>0$ for all $x\in\Omega$. Since $\varphi\in C(\overline{\Omega})$, we deduce that there exists $c>0$ such that $\varphi(x)\ge c$ for all $x\in\overline{\Omega}$.

    The $G$-invariance and strict positivity of $\varphi$ also imply the simplicity of the eigenvalue $\lambda_0$. 
    Indeed, suppose for contradiction that there exists an eigenfunction $\tilde\varphi$ corresponding to $\lambda_0$ and orthogonal to $\varphi$. 
    The arguments from the preceding paragraphs also apply to $\tilde\varphi$, so in particular $\tilde \varphi \ge 0$ or $\tilde \varphi \le 0$. 
    This contradicts the orthogonality of $\varphi$ and $\tilde \varphi$. 
    
    Finally, let $\mu\in (0, \lambda_0)$ be another positive eigenvalue of $-L_B$ with a normalised eigenfunction $u$.  
    Then we can write $u = u_0 + u_1$ where $u_0 \in F_G \cap H^1(\Omega)$ and $u_1 \in F_G^\perp \cap H^1(\Omega)$, see Proposition~\ref{prop:symmetry}\ref{prop:symmetry:itm:H1-directsum} in the appendix. 
    Since $L_B$ is $G$-equivariant, it follows from Proposition~\ref{prop:symmetry}\ref{prop:symmetry:itm:op} that $u_0 \in \dom(L_B)$ and $-L_B u_0 = \mu u_0$. 
    Let us show that $u_0 \neq 0$. 
    To this end, note that
    \[
        0 
        > 
        -\mu 
        = 
        \mf{a}_B[u,u]
        = 
        \mf{a}_B[u_0,u_0] + \mf{a}_B[u_1,u_1]
        \ge 
        \mf{a}_B[u_0,u_0]
        ;
    \]
    here, the first equality follows from $-L_Bu = \mu u$, 
    the second equality is due to Proposition~\ref{prop:symmetry}\ref{prop:symmetry:itm:form}, 
    and the inequality at the end holds according to~\eqref{thm:symm-eig:eq:aB-lower-bound} since $u_1 \in F_G^\perp \cap H^1(\Omega)$.
    Thus, $u_0 \neq 0$ as claimed. 
    
    The same argument that we used to show that $\varphi(x) > 0$ for all $x \in \overline{\Omega}$ can also be applied to $u_0$ since $\mu > 0$. 
    Hence, $u_0 \ge 0$ or $u_0 \le 0$. 
    But the self-adjointness of $L_B$ and $\mu \ne \lambda_0$ imply that $u_0$ is orthogonal to $\varphi$, which is a contradiction.
    
    \ref{thm:ball:itm:ev-pos}
    The semigroup $(e^{-tL_B})_{t\ge 0}$ is self-adjoint and real by Proposition~\ref{prop:basic}\ref{prop:basic:itm:self-adjoint} and~\ref{prop:basic:itm:real} respectively. Due to the assumptions that $B$ is self-adjoint on $L^2(\partial\Omega)$ and leaves $L^\infty(\partial\Omega)$ invariant, we deduce from Lemma~\ref{lem:op-L_1-L_infty}\ref{lem:op-L_1-L_infty:itm:adjoint-2} that $B$ acts boundedly on $L^1$ and $L^\infty$. Hence Theorem~\ref{thm:ultra} implies that $e^{-tL_B}(L^2(\Omega))\subseteq L^\infty(\Omega)$ for all $t>0$.
	
	According to part~\ref{thm:ball:itm:spectrum}, we can choose an eigenfunction $\psi$ associated to $\rs(-L_B)$ satisfying $\psi \ge \mathbf{1}$. Since $\rs(-L_B)$ is a simple eigenvalue and $-L_B$ is self-adjoint, the desired conclusion follows from Proposition~\ref{prop:unif-ev-pos}.
\end{proof}

\begin{remark}
	In the setting of Theorem~\ref{thm:ball}, we can choose an eigenfunction $\varphi \in L^\infty(\Omega)$ associated to the eigenvalue $\rs(-L_B)$ such that $\varphi\ge \mathbf{1}$. Combined with the ultracontractivity property $e^{-tL_B}(L^2(\Omega)) \subseteq L^\infty(\Omega)$, this shows that for each $t>0$ there exists a constant $c_t>0$ such that $\modulus{e^{-tL_B}f} \le c_t \varphi$ for all $f\in L^2(\Omega)$. In the terminology introduced by Davies and Simon~\cite{DS84}, this implies that the semigroup $(e^{-tL_B})_{t\ge 0}$ is \emph{intrinsically ultracontractive}. This property is well-known for a large variety of second-order elliptic operators with \emph{local} boundary conditions.
\end{remark}

\begin{example}
\label{exa:bose-model}
    We revisit the model of the Bose condensation example in~\cite[Theorem 6.13]{DGK2} and prove uniform eventual positivity of the semigroup under a different set of assumptions on the boundary operator $B$.
    
    In Setting~\ref{setting:main}, let $\Omega = \{ x\in\RR^2 ~\colon \modulus{x} < 1\}$ be the open unit disk and suppose $A(x) = \id_{\RR^2}$ for all $x\in\Omega$ (i.e.\ we consider the Laplacian). 
    One can identify the boundary $\partial \Omega$ with the complex unit circle, so it becomes a multiplicative group.
    Let $q \in L^1(\partial\Omega)$ be a real-valued function such that its Fourier coefficients, defined by
    \begin{equation*}
        q_k \coloneqq \int_{\partial \Omega} q(z) z^{-1} \, dz \qquad\forall\, k\in\ZZ,
    \end{equation*}
    are real for all $k\in\ZZ$, and $q_0 < 0$. 
    Define the convolution operator $B \in \mathcal{L}(L^2(\partial \Omega))$ by
    \begin{equation*}
	   (Bf)(z) \coloneqq (q*f)(z) \coloneqq \int_{\partial \Omega} q(w)f(w^{-1} z) \, dw
    \end{equation*}
    for all $f \in L^2(\partial \Omega)$ and for  $z \in \partial \Omega$. 
    Consider the special orthogonal group $G = SO(2)$, which can be identified with the multiplicative group $\partial \Omega$ that acts on itself.
    Then the following assertions hold:
    \begin{enumerate}[(i)]
        \item\label{exa:bose-model:itm:setting} 
        For the choice of coefficient matrix $A$, the operator $B$, and the group $G$ given above, all the conditions of Setting~\ref{setting:symmetry} are satisfied. Moreover we have $\braket{B\mathbf{1}_{\partial\Omega}, \mathbf{1}_{\partial\Omega}}_{\partial\Omega} < 0$.

        \item\label{exa:bose-model:itm:ev-pos} 
        There exists a constant $\beta > 0$ (independent of $q$) such that if $\|q\|_1 \le \beta$, then there exist $t_0 \ge 0$ and a constant $\delta > 0$ such that
        \[
            e^{-tL_B}f \ge \delta \left( \int_\Omega f \,dx \right)\mathbf{1} \qquad\forall\, t\ge t_0
        \]
        for all $0\le f\in L^2(\Omega)$. 
        So in particular, the semigroup $(e^{-tL_B})_{t \ge 0}$ is uniformly eventually positive if $\|q\|_1 \le \beta$.
    \end{enumerate}

    \begin{proof}
        \ref{exa:bose-model:itm:setting}
        Clearly $\Omega$ is a $G$-invariant Lipschitz domain and the identity matrix satisfies~\ref{setting:symmetry:itm:coefficients}. 
        Moreover, $B$ satisfies condition~\ref{setting:symmetry:itm:B2}. 
        Indeed, $B$ is self-adjoint since the Fourier coefficients of $q$ are real. 
        The $G$-equivariance of $B$ is a classical property of convolutions and can be seen as follows: 
        for $y \in G \simeq \partial \Omega$ and the translation operator $S_y \in \mathcal{L}(L^2(\partial \Omega))$ introduced in Section~\ref{sec:invariance-group-actions} one has
        \[
            (BS_y) f(z) 
            = 
            \int_{\partial \Omega} q(w)(S_y f)(w^{-1}z)\,dw 
            = 
            \int_{\partial \Omega} q(w)f(w^{-1}y^{-1}z)\,dw 
            = 
            (S_y Bf)(z)
        \]
        for all $f \in L^2(\partial \Omega)$ and for $z \in \partial \Omega$.
        Thus, all assumptions of Setting~\ref{setting:symmetry} are satisfied. 
        Finally, observe that
        \[
            \braket{B\mathbf{1}_{\partial\Omega}, \mathbf{1}_{\partial\Omega}}_{\partial\Omega} = 2\pi\int_{-\pi}^\pi q(x)\,dx = 2\pi q_0 < 0.
        \]

        \ref{exa:bose-model:itm:ev-pos}
        We check the conditions of Theorem~\ref{thm:ball}.
        Since $q$ is real-valued, it follows that $B$ is a real operator on $L^2(\partial\Omega)$.
        Obviously $B$ leaves $L^\infty(\partial\Omega)$ invariant, and $G$ acts transitively on $\partial\Omega = \mathbb{S}^1$. Together with the results of part~\ref{exa:bose-model:itm:setting}, we have shown that all the conditions of Theorem~\ref{thm:ball} are indeed satisfied. 
        Since Young's convolution inequality yields $\|B\|_{L^2(\partial\Omega) \to L^2(\partial\Omega)} \le \|q\|_1$, the assertion now follows from Theorem~\ref{thm:ball}\ref{thm:ball:itm:ev-pos}.
    \end{proof}
\end{example}

We note in passing that, at the end of the preceding proof, one actually has the equality $\|B\|_{L^2(\partial\Omega) \to L^2(\partial\Omega)} = \|q\|_1$, as one can see by applying $B$ to constant functions.

We do not investigate whether the smallness condition~\eqref{eq:small-B} in Theorem~\ref{thm:ball} is optimal. However, the following computations show that Theorem~\ref{thm:ball} is not true without some upper bound on the norm of $B$.

\begin{example}
	In this example, we consider the $1$-dimensional ball $\Omega = (-1,1)$. The group $G=O(1)$ of order $2$ acts on $\RR$ by reflection, and clearly $\Omega$ is $G$-invariant. The space $L^2(\partial\Omega)$ may be identified with $\CC^2$, and $G$ acts on $\CC^2$ by permuting the coordinates. Define the boundary operator
	\begin{equation*}
		B = b\begin{pmatrix} -2 & 1 \\ 1 & -2 \end{pmatrix}
	\end{equation*}
	with a parameter $b>0$. One easily verifies that $B$ satisfies~\ref{setting:symmetry:itm:B2} in Setting~\ref{setting:symmetry} --- in particular it commutes with the permutation matrix $\left(\begin{smallmatrix} 0 & 1 \\ 1 & 0 \end{smallmatrix}\right)$, so it is $G$-equivariant. Our objective is to show explicitly that if $b>0$ is sufficiently large, then the conclusions of Theorem~\ref{thm:ball} fail to hold.
	
	For $\lambda>0$, we solve the eigenvalue problem
	\begin{equation*}
		u'' = \lambda u \quad\text{in }\Omega, \qquad \begin{pmatrix} -u'(-1) \\ u'(1) \end{pmatrix} + B\begin{pmatrix} u(-1) \\ u(1) \end{pmatrix} = \begin{pmatrix} 0 \\ 0 \end{pmatrix}.
	\end{equation*}
	A general solution to the eigenvalue equation is given by
	\begin{equation*}
		u(x) = c_1 \cosh(\mu x) + c_2 \sinh(\mu x), \qquad \mu\coloneqq\sqrt{\lambda}.
	\end{equation*}
	It is helpful to observe that
	\begin{align*}
		\begin{pmatrix} -u'(-1) \\ u'(1) \end{pmatrix} &= \begin{pmatrix} \mu\sinh\mu & -\mu\cosh\mu \\ \mu\sinh\mu & \mu\cosh\mu \end{pmatrix} \begin{pmatrix} c_1 \\ c_2 \end{pmatrix} \quad\text{and} \\
		\begin{pmatrix} u(-1) \\ u(1) \end{pmatrix} &= \begin{pmatrix} \cosh\mu & -\sinh\mu \\ \cosh\mu & \sinh\mu \end{pmatrix} \begin{pmatrix} c_1 \\ c_2 \end{pmatrix}.
	\end{align*}
	After some simple computations, the boundary conditions can be rewritten as
	\begin{equation}
	\label{eq:BCs-mu}
		\begin{pmatrix} \mu\sinh\mu - b\cosh\mu & 3b\sinh\mu - \mu\cosh\mu \\ \mu\sinh\mu - b\cosh\mu & \mu\cosh\mu - 3b\sinh\mu \end{pmatrix} \begin{pmatrix} c_1 \\ c_2 \end{pmatrix} = \begin{pmatrix} 0 \\  0 \end{pmatrix}.
	\end{equation}
	Therefore, the condition for $\lambda$ to be an eigenvalue of $-L_B$ is
	\begin{equation*}
	\begin{aligned}
		\det &\begin{pmatrix} \mu\sinh\mu - b\cosh\mu & 3b\sinh\mu - \mu\cosh\mu \\ \mu\sinh\mu - b\cosh\mu & \mu\cosh\mu - 3b\sinh\mu \end{pmatrix} \\
		&\qquad\qquad = 2(\mu\sinh\mu - b\cosh\mu)(\mu\cosh\mu - 3b\sinh\mu) = 0.
	\end{aligned}
	\end{equation*}
	Consider the two functions given by
	\begin{equation}
		f_1(\mu) = \mu\tanh\mu \quad\text{and}\quad f_2(\mu) = \frac{1}{3}\mu\coth\mu
	\end{equation}
	for $\mu \ge 0$, which are plotted in Figure~\ref{fig:mu-plots}.
        Then $\lambda = \mu^2$ is an eigenvalue of $-L_B$ if and only if $f_1(\mu) = b$ or $f_2(\mu) = b$. 
	
	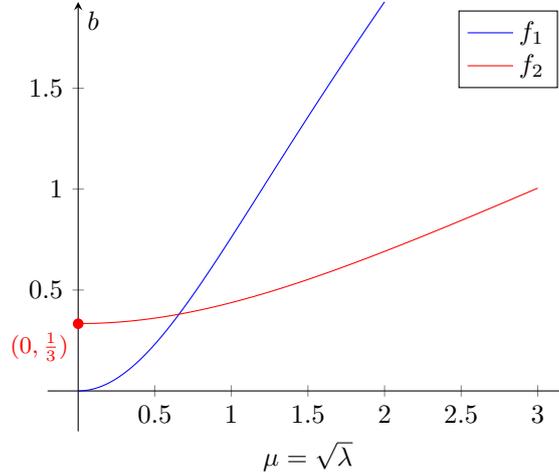
\begin{figure}[h]
		\caption{The graphs of $f_1, f_2$ encode the positive eigenvalues for a given $b>0$.}
		\centering
		\begin{tikzpicture}
			\label{fig:mu-plots}
			\begin{axis}[axis lines = middle, clip=false, 
				xmin=-0.2, xmax=3.2, ymin=-0.2, xlabel = {$\mu=\sqrt{\lambda}$}, xlabel near ticks, ylabel = $b$]
				
				\addplot[domain = 0:2, samples=80, color=blue, dashed]
				{x*tanh(x)};
				\addlegendentry{$f_1$};
				
				\addplot[domain = 0.02:3, samples=80, color=red]
				{x/(3*tanh(x))};
				\addlegendentry{$f_2$};
				
				\fill[red] (0,0.333) circle (2pt) node[below left] {\small $(0, \frac{1}{3})$};
			\end{axis}
		\end{tikzpicture}
	\end{figure}
	
	Observe that $f_2(0) = \frac{1}{3}$, and when $b<\frac{1}{3}$, there is precisely one positive eigenvalue; namely $\lambda=\mu^2$, where $\mu$ is the unique positive solution to the equation $b=f_1(\mu) = \mu\tanh\mu$. In that case, the first column of the matrix in equation~\eqref{eq:BCs-mu} consists of zeroes, and thus we can take $c_2=0$ and $c_1>0$ arbitrary to obtain the eigenfunction
	\begin{equation*}
		\varphi(x) = c_1 \cosh(\mu x), \quad \mu = f_1^{-1}(b).
	\end{equation*}
	Evidently $\varphi$ is strictly positive and $G$-invariant on $[-1,1]$.
	
	When $b=\frac{1}{3}$, a second solution $\mu=0$ (yielding the eigenvalue $\lambda=0$) appears. By simple calculations, one checks readily that the $0$-eigenspace is spanned by the function $\psi(x)=x$, which is notably not positive on $[-1,1]$. The curves $f_1, f_2$ cross at the value
	\begin{equation*}
		\mu^* = \tanh^{-1}\left(\frac{1}{\sqrt{3}}\right) \approx 0.658479,
	\end{equation*}
	in which case the corresponding eigenspace is two-dimensional and spanned by $\{\cosh(\mu^* \,\cdot\,), \sinh(\mu^* \,\cdot\,)\}$. As $b$ increases further, the double eigenvalue splits. However, the leading eigenvalue now arises from $f_2$, and one can verify using~\eqref{eq:BCs-mu} that the corresponding eigenspace is spanned by $\sinh(\mu \,\cdot\,)$. Thus the positivity of the leading eigenfunction is lost, and it is the smaller positive eigenvalue arising from $f_1$ that yields a positive eigenfunction.
\end{example}

\subsection*{Open problem}

In Setting~\ref{setting:main}, consider the special case of the Laplacian, which corresponds to $A(x) = \id_{\RR^d}$ for all $x\in\Omega$. As far as spectral theory is concerned, the results of this article are complementary to the analysis in~\cite{DGK2}, in the sense that we have developed a general theory under the condition $\rs(-L_B)\ge 0$, whereas the examples considered in Theorems~6.11 and~6.13 in the aforementioned paper satisfy $\rs(-L_B)<0$. For specific choices of $B$, uniform eventual positivity of the semigroup $(e^{-tL_B})_{t\ge 0}$ was proved by showing an analogous positivity property of the resolvent operator 
$(\lambda + L_B)^{-1}$ as $\lambda \downarrow 0$, and this was achieved by explicit computations. 
It would be desirable to analyse the case of a negative spectral bound in a systematic way, but at the time of writing, it is not clear to the authors how to achieve this.

\subsection*{Acknowledgements}

This research was funded by the Deutsche Forschungsgemeinschaft (DFG, German Research Foundation) – Project No.\ 515394002.

We thank the referees for many comments and suggestions that greatly improved the manuscript.

\appendix

\section{Group actions and their invariant functions}
\label{app:symmetric}

In Section~\ref{sec:symm}, the action of a closed subgroup $G \subseteq O(d)$ on a $G$-invariant set $\Omega \subseteq \RR^d$ occurs. 
The following proposition shows several properties of such group actions that are used in the proof of Theorem~\ref{thm:symm-eig}.
Recall that on every compact group $G$, there exists a unique normalised Haar measure $\mu$. 
It is both left and right invariant, see e.g.~\cite[Corollary 2.28]{Folland}.

\begin{proposition}
    \label{prop:symmetry}
    Let $\emptyset \not= \Omega \subseteq \RR^d$ be open, let $G \subseteq O(d)$ be a closed subgroup that leaves $\Omega$ invariant, and let $\mu$ be the normalised Haar measure on $G$. 
    Let $F_G \subseteq L^2(\Omega) \coloneqq L^2(\Omega)$ denote the subspace of $G$-invariant functions (see~\eqref{eq:symm-subspace}).
    
    \begin{enumerate}[\upshape (i)]
        \item\label{prop:symmetry:itm:cont-group-action} 
        The group actions
        \begin{align*}
            G \times L^2(\Omega) &\to L^2(\Omega), \quad 
            (g, f) \mapsto T_g f = f(g^{-1}\,\cdot\,) \\
            G \times H^1(\Omega) &\to H^1(\Omega), \quad 
            (g, u) \mapsto T_g u = u(g^{-1}\,\cdot\,)
        \end{align*}
        are continuous maps.
        
        \item\label{prop:symmetry:itm:symmetric-proj-L2}  
        Let $P \in \mathcal{L}(L^2(\Omega))$ be the orthogonal projection onto $F_G$. 
        Then 
        \begin{equation}
            \label{prop:symmetry:eq:symmetric-proj}
            Pf = \int_G T_g f \,d\mu(g) \qquad\forall\, f\in L^2(\Omega),
        \end{equation}
        where the integral exists as a Bochner integral with values in $L^2(\Omega)$.

        \item\label{prop:symmetry:itm:symmetric-proj-H1} 
        If $f \in H^1(\Omega)$, then the integral in~\eqref{prop:symmetry:eq:symmetric-proj} exists as a Bochner integral in $H^1(\Omega)$, and hence $Pf \in H^1(\Omega)$.

        \item\label{prop:symmetry:itm:form}
        If $\mf{c} \colon H^1(\Omega)\times H^1(\Omega) \to \CC$ is a bounded and $G$-invariant sesquilinear form, then
        \[
            \mf{c}[u-Pu,Pu] = 0 \qquad\forall\,u\in H^1(\Omega).
        \]

        \item\label{prop:symmetry:itm:op}
        If $C \colon \dom(C) \subseteq L^2(\Omega) \to L^2(\Omega)$ is a closed and $G$-equivariant operator, then $Pf \in \dom(C)$ and $CPf = PCf$ for all $f \in \dom(C)$.

        \item\label{prop:symmetry:itm:H1-directsum}
        One has the direct sum decomposition
        \[
            H^1(\Omega) = F_G\cap H^1(\Omega) \oplus_{H^1} F_G^\perp \cap H^1(\Omega),
        \]
        where the summands are orthogonal with respect to the $H^1$ inner product.
    \end{enumerate}

    From now on, assume in addition that $\Omega$ is a bounded Lipschitz domain, and let $\gamma \colon H^1(\Omega) \to L^2(\partial\Omega)$ denote the trace operator.
    
    \begin{enumerate}[\upshape (i), resume]
        \item\label{prop:symmetry:itm:bdry-intertwine} 
        One has $\gamma (T_g u) = S_g \gamma(u)$ for every $u \in H^1(\Omega)$.
        
        \item\label{prop:symmetry:itm:symmetric-proj-bdry}
        The orthogonal projection $Q \in \mathcal{L}(L^2(\partial\Omega))$ onto the subspace of $G$-invariant functions in $L^2(\partial\Omega)$ is given by
        \[
            Qf = \int_G S_g f \,d\mu(g) \qquad\forall\, f\in L^2(\partial\Omega).
        \]
        Hence, $\gamma (P u) = Q \gamma(u)$ for every $u \in H^1(\Omega)$, where $P \in \mathcal{L}(L^2(\Omega))$ is the projection from~\ref{prop:symmetry:itm:symmetric-proj-L2}.
        
        \item\label{prop:symmetry:itm:trace} 
        The spaces $\gamma(F_G\cap H^1(\Omega))$ and $\gamma(F_G^\perp \cap H^1(\Omega))$ are orthogonal in $L^2(\partial \Omega)$. 
        In particular, 
        $
            \braket{\one_{\partial \Omega}, \gamma(v)}_{\partial \Omega} = 0
        $
        for every $v \in F_G^\perp \cap H^1(\Omega)$.
    \end{enumerate}
\end{proposition}

\begin{proof}
    \ref{prop:symmetry:itm:cont-group-action}
    Let $e$ denote the unit element of $G$ (i.e.\ the identity matrix in $\RR^{d\times d}$). It suffices to prove continuity at $(e,f_0) \in G \times L^2(\Omega)$, respectively at $(e,u_0) \in G \times H^1(\Omega)$, for arbitrary $f_0 \in L^2(\Omega)$ and $u_0 \in H^1(\Omega)$.
    
    Firstly, suppose $f\in C_c(\Omega)$. Then by the uniform continuity of $f$ and continuity of the inversion $g \mapsto g^{-1}$, it follows that $T_g f \to f$ in $L^2(\Omega)$ as $g \to e$. Since $C_c(\Omega)$ is dense in $L^2(\Omega)$ and $\|T_g\|_{L^2 \to L^2} = 1$ for all $g\in G$, the convergence then holds for all $f\in L^2(\Omega)$.
    Now let $f_0, f_1 \in L^2(\Omega)$ be arbitrary. The desired continuity of the group action on $L^2(\Omega)$ then follows from the inequality
    \[
        \|T_g f_1 - f_0\|_2 \le \|T_g (f_1 - f_0)\|_2 + \|T_g f_0 - f_0\|_2 .
    \]

    To prove the continuity of the action on $H^1(\Omega)$, we recall that each operator $T_g$ leaves $H^1(\Omega)$ invariant and $\|T_g\|_{H^1 \to H^1} = 1$; both properties follow from the chain rule for Sobolev functions, as mentioned at the beginning of Section~\ref{sec:invariance-group-actions}.
    For all $u_0, u_1 \in H^1(\Omega)$, we obtain
    \begin{align*}
        \| \nabla T_g u_1 - \nabla u_0\|_2 &\le \| \nabla T_g(u_1 - u_0)\|_2 + \| \nabla T_g u_0 - \nabla u_0\| \\
        &\le \|u_1 - u_0\|_{H^1} + \| g (\nabla u_0 \circ g^{-1}) - \nabla u_0\|_2,
    \end{align*}
    and the desired continuity on $H^1(\Omega)$ now follows easily from the continuity of the group action on $L^2(\Omega)$ and the above inequality.

    \ref{prop:symmetry:itm:symmetric-proj-L2}
    For each $f \in L^2(\Omega)$ the map $G \ni g \mapsto T_g f \in L^2(\Omega)$ is continuous by part~\ref{prop:symmetry:itm:cont-group-action} and is thus Bochner integrable.
    One clearly has $F_G \subseteq P(L^2(\Omega))$ since $\mu(G) = 1$. 
    The translation invariance of the Haar measure gives the converse inclusion $P(L^2(\Omega)) \subseteq F_G$ and that $P$ is a projection. Finally, $P$ is contractive and thus orthogonal since every operator $T_g$ is contractive.

    \ref{prop:symmetry:itm:symmetric-proj-H1}
    Let $f \in H^1(\Omega)$.
    Part~\ref{prop:symmetry:itm:cont-group-action} shows that the integrand in~\eqref{prop:symmetry:eq:symmetric-proj} is continuous and thus Bochner integrable with values in $H^1(\Omega)$. 
    As $H^1(\Omega)$ embeds continuously into $L^2(\Omega)$, the value of the Bochner integral in both spaces coincides, so indeed $Pf \in H^1(\Omega)$.
    
    \ref{prop:symmetry:itm:form}
    Let $u\in H^1(\Omega)$. The $G$-invariance of $\mf{c}$ and part~\ref{prop:symmetry:itm:symmetric-proj-H1} together imply
    \begin{align*}
        \mf{c}[u-Pu,Pu] = \int_G \mf{c}[u-Pu,Pu] \,d\mu(g) &= \int_G \mf{c}[T_g(u-Pu),Pu] \,d\mu(g) \\
        &= \mf{c}[P(u-Pu), Pu] = 0;
    \end{align*}
    in the second equality we have used that $T_g Pu = Pu$ and in third equality we have used that $\mf{c}[\,\cdot\,,v]$ is a bounded linear functional on $H^1(\Omega)$ for each fixed $v\in H^1(\Omega)$.

    \ref{prop:symmetry:itm:op} 
    Let $f \in \dom(C)$.
    The continuity of the map $G \to L^2(\Omega)$, $g \mapsto T_g h$ for every $h \in L^2(\Omega)$ and the $G$-equivariance of $C$ imply that $G \to \dom(C)$, $h \mapsto T_g f$ is continuous with respect to the graph norm on $\dom(C)$. 
    Since $C$ is closed, $\dom(C)$ is a Banach space with respect to the graph norm, 
    so the integral~\eqref{prop:symmetry:eq:symmetric-proj} in $L^2(\Omega)$ that represents $Pf$ also converges as a Bochner integral in $\dom(C)$. 
    As the embedding of $\dom(C)$ into $L^2(\Omega)$ is continuous, the integrals in both spaces coincide, so indeed $Pf \in \dom(C)$.  
    The formula $CPf = PCf$ now follows from the $G$-equivariance of $C$ and the fact that $C$ is continuous from $\dom(C)$ (equipped with the graph norm) to $L^2(\Omega)$.
    
    \ref{prop:symmetry:itm:H1-directsum}
    We first check that the standard $H^1(\Omega)$ inner product
    \[
        (u,v) \mapsto \int_\Omega u\overline{v}\,dx + \int_\Omega \nabla u \cdot \overline{\nabla v}\,dx
    \]
    is a $G$-equivariant sesquilinear form. This property is clear for the first summand in the above inner product. 
    For the second summand it follows from the chain rule for Sobolev functions, which gives $\nabla (u\circ g^{-1}) = g (\nabla u \circ g^{-1})$ for all $g\in G$ and $u\in H^1(\Omega)$ (cf.\ the proof of part~\ref{prop:symmetry:itm:cont-group-action}). 
    Thus, one can apply~\ref{prop:symmetry:itm:form} to the $H^1(\Omega)$ inner product to conclude that the spaces $F_G\cap H^1(\Omega)$ and $F_G^\perp\cap H^1(\Omega)$ are orthogonal in $H^1(\Omega)$. 
    Moreover, their sum equals $H^1(\Omega)$ since this space is invariant under $P$ by part~\ref{prop:symmetry:itm:symmetric-proj-H1}.
    \smallskip 

    For the rest of the proof, assume that $\Omega$ is a bounded Lipschitz domain.
    \smallskip 
    
    \ref{prop:symmetry:itm:bdry-intertwine}
    Let $g \in G$. 
    For $u\in C(\overline{\Omega})\cap H^1(\Omega)$ the equality $\gamma (T_g u) = S_g \gamma(u)$ is clear. 
    Since $\Omega$ is a Lipschitz domain, the equality for $u\in H^1(\Omega)$ follows by density of $C(\overline{\Omega})\cap H^1(\Omega)$ in the $H^1$ norm (see for instance~\cite[Theorem 11.35]{Leoni}) and the continuity of the operator $T_g$ from $H^1(\Omega)$ to $H^1(\Omega)$.
    
    \ref{prop:symmetry:itm:symmetric-proj-bdry}
    The proof of the formula for $Q$ is analogous to the proof of~\eqref{prop:symmetry:eq:symmetric-proj} in~\ref{prop:symmetry:itm:symmetric-proj-L2}.
    For all $u \in H^1(\Omega)$ one thus obtains
    \[
        \gamma (Pu) = \gamma \left( \int_G T_g u \,d\mu(g) \right) = \int_G \gamma(T_g u) \,d\mu(g) = \int_G S_g\gamma(u) \,d\mu(g) = Q\gamma(u)
    \]
    where we used~\ref{prop:symmetry:itm:symmetric-proj-L2} for the first equality,~\ref{prop:symmetry:itm:symmetric-proj-H1} for the second equality, and~\ref{prop:symmetry:itm:bdry-intertwine} for the third equality.
    
    \ref{prop:symmetry:itm:trace}
    Let $P \in \mathcal{L}(L^2(\Omega))$ and $Q \in \mathcal{L}(L^2(\partial \Omega))$ denote the orthogonal projections from~\ref{prop:symmetry:itm:symmetric-proj-L2} and~\ref{prop:symmetry:itm:symmetric-proj-bdry} respectively, and let $u\in F_G \cap H^1(\Omega)$ and $v\in F_G^\perp \cap H^1(\Omega)$. 
    It follows from~\ref{prop:symmetry:itm:symmetric-proj-bdry} that $\gamma(u) = \gamma(Pu) = Q\gamma(u)$ and $Q\gamma(v) = \gamma(Pv) = 0$, so
    \begin{align*}
        \braket{\gamma(u), \gamma(v)}_{\partial\Omega} 
        = 
        \braket{Q\gamma(u), \gamma(v)}_{\partial \Omega}
        = 
        \braket{\gamma(u), Q\gamma(v)}_{\partial \Omega}
        =
        0
        ,
    \end{align*}
    where we used the self-adjointness of $Q$ for the second equality.

    The equality $\braket{\one_{\partial \Omega}, \gamma(v)}_{\partial \Omega} = 0$ for all $v \in F_G^\perp \cap H^1(\Omega)$ now follows since $\one_{\partial \Omega} = \gamma(\one) \in \gamma(F_G \cap H^1(\Omega))$.
\end{proof}

\bibliographystyle{plain}
\bibliography{references}

\end{document}